\documentclass[review]{elsarticle}

\usepackage{hyperref}
\journal{Chaos}

\usepackage{amsmath}
\usepackage{amsthm}
\usepackage{tikz-cd}
\usepackage{circuitikz}
\usepackage{algorithmic}
\usepackage{algorithm}
\usepackage{dsfont}
\usepackage{color}
\usepackage{amssymb}
\usepackage[english]{babel}
\usepackage{graphicx,color}
\usepackage{extarrows}
\usepackage{mathrsfs}
\usepackage{booktabs}
\usepackage{multirow}
\usepackage{upgreek}
\usepackage{threeparttable}
\usepackage{textcomp}
\usepackage[justification=raggedright]{caption}
\usepackage[misc]{ifsym}
\usepackage{array}
\usepackage{graphicx}
\usepackage{subfloat}
\usepackage{float}
\usepackage{subfigure}

\newtheorem{theorem}{Theorem}
\usepackage[T1]{fontenc}

\makeatletter
\@addtoreset{equation}{section}
\makeatother

\graphicspath{ {./figures/} }

\begin{document}
        \begin{frontmatter}

            \title{Causal Discovery in Symmetric Dynamic Systems with Convergent Cross Mapping}

            \author[inst1]{Yiting Duan\corref{mycorrespondingauthor}}
            \cortext[mycorrespondingauthor]{Corresponding author}
            \ead{Y.Duan@westernsydney.edu.au}
		\affiliation[inst1]{organization={School of Computer, Data and Mathematical Sciences},
			addressline={Western Sydney University}, 
			city={Parramatta, NSW},
			postcode={2150},
			country={Australia}}

            \author[inst1]{Yi Guo}
		\ead{y.guo@westernsydney.edu.au}
		
		\author[inst2]{Jack Yang}
            \ead{jiey@uow.edu.au}
		\affiliation[inst2]{organization={School of Computing and Information Technology},
			addressline={University of Wollongong}, 
			city={Wollongong},
			postcode={2500},
			country={Australia}}

		\author[inst3]{Ming Yin}
            \ead{m.yin@scnu.edu.cn}
            \affiliation[inst3]{organization = {School of Semiconductor Science and Technology},
                addressline = {South China Normal University},
                city = {Guangzhou},
                postcode = {510631},
                country = {China}}

        \begin{abstract}
            This paper systematically discusses how the inherent properties of chaotic attractors influence the results of discovering causality from time series using convergent cross mapping, particularly how convergent cross mapping misleads bidirectional causality as unidirectional when the chaotic attractor exhibits symmetry. We propose a novel method based on the k-means clustering method to address the challenges when the chaotic attractor exhibits two-fold rotation symmetry. This method is demonstrated to recover the symmetry of the latent chaotic attractor and discover the correct causality between time series without introducing information from other variables. We validate the accuracy of this method using time series derived from low-dimension and high-dimensional chaotic symmetric attractors for which convergent cross mapping may conclude erroneous results. 
            
        \end{abstract}

            
        \end{frontmatter}
        
        \textbf{Discovering nonlinear causal interactions in complex chaotic systems is essential for uncovering the underlying physical phenomena. Convergent cross mapping is the crucial tool for this causality discovery which is based on embedding the time series into a shadow manifold that is topologically equivariant to the original attractor through delay embedding. However, despite the problem of selecting embedding parameters, the accuracy of the obtained causality is also highly related to the intrinsic properties of the original attractor, especially when the chaotic attractor exhibits rotation symmetry, such as Lorenz63 and Burke $\&$ Shaw. Under this circumstance, the shadow manifold reconstructed from the invariant variable may lose the original symmetry, causing the convergent cross mapping to mislead bidirectional causality as unidirectional. In this work, we illustrate this as the reconstructing mapping mods out the symmetry, thereby no longer being a one-to-one mapping. We propose a new framework based on k-means clustering to solve this problem without introducing any information from another variable. The efficiency of this method is demonstrated using various time series derived from both low-dimensional and high-dimensional chaotic symmetric attractors.}  
        
    \section{Introduction}

    Identifying the causation between signals from complex systems has garnered significant interest and is essential in ecology, epidemiology, and climatology, where causal relations are inferred from observational signals. Granger causality, as outlined by \cite{granger1969investigating}, is a prevalent method for detecting causality in signals based on the principle of predictability and has been widely used to analyze economic growth and neuroscience for effective brain connectivity \cite{kitole2024navigating, ross2024causation}. However, the effectiveness of Granger causality is contingent on the separability of the system, making it inadequate for signals from nonlinear or deterministic physical systems, particularly systems with attractors characterized by non-separability and moderate coupling.

    Addressing the limitations of Ganger causality, Sugihara \cite{sugihara2012detecting} proposed the Convergent Cross Mapping(CCM) method. This method leverages the principles of chaotic dynamics to detect causality between signals from nonlinear, nonseparable, and moderately coupled dynamic systems. The key idea of CCM for detecting causality between two signals is determining whether the historical states of one signal can reliably predict the states of another by checking the existence of the mapping between two attractors reconstructed from two signals. Recently, CCM has been successfully applied in several fields, including climatology \cite{qing2023soil}, neuroscience \cite{avvaru2023effective}, and ecology \cite{gao2023causal}.

    The application of CCM relies heavily on the embedding technique, mainly the delay-coordinate mapping, which serves as a diffeomorphism between the reconstructed manifold, also termed the shadow manifold $\mathcal{M}$, and the original attractor manifold. Ideally, this diffeomorphism ensures that the dynamics of the original system are faithfully captured in the shadow manifold. However, the delay-coordinate mapping can not reach a diffeomorphism in practical scenarios. This discrepancy can arise due to various external factors, including noise, missing values, and improper selection of reconstruction parameters, e.g., the inappropriate lags for the delay-coordinate mapping. Furthermore, we point out that even in cases where the signals are clean and sufficiently lengthy, intrinsic properties of the chaotic system, such as the observability of the variables and inherent symmetries, may also influence the embedding quality. Despite the significance of these factors, few studies have thoroughly examined the impact of such inherent factors on CCM results.

    Moreover, special attention should be given to cases when signals from the attractor exhibit rotation symmetry since the shadow manifold reconstructed from the signals of the invariant variable of this attractor may lose the symmetry such that the delay-coordinate mapping is no longer one-to-one. For example, the shadow manifold reconstructed from the $z$-coordinate projection of the Lorenz63 system does not exhibit the same butterfly attractor as the original ones, as noted in counterexamples provided in \cite{yuan2020data, butler2023causal} and the supplementary material of \cite{sugihara2012detecting}. Although such time series arising from rotation symmetric attractors are abundant in real scenarios, such as output signals from chaotic circuits, the previous works have yet offer a practical solution to address this issue.

    In this work, we make the following contributions:
    \begin{enumerate}
        \item We systematically discuss how the inherent properties of chaotic systems, particularly those with symmetric characteristics, influence the results of CCM and derive a general conclusion.   
        \item We propose a new framework based on the k-means clustering method to address the challenges of two-fold rotation symmetry.
        \item We validate the correctness and efficiency of our method via several examples.     
    \end{enumerate}
    
    The rest of this article is organized as follows. Section \ref{SectionTwo} introduces CCM and discusses related work in this field. Section \ref{SectionThree} describes the challenge of CCM for several counterexamples. In section \ref{SectionFour}, we explain the reason by systematically discussing the influence of symmetry and propose our framework. Section \ref{SectionFive} assesses the proposed method in several cases. Finally, conclusions are discussed in Section \ref{Conclusion}.

    \section{Background}\label{SectionTwo}
    
    In this section, we provide a self-contained introduction to the CCM methodology, specifically focusing on its ability to identify causality between two variables \(X\) and \(Y\), generated from two deterministic dynamic systems. 
    
    \subsection{Attractor Reconstruction}

    In practical applications, signals from dynamic systems with attractors are often shown as time series of real numbers generated from observations. The goal is to reconstruct the state space of the unknown attractor for analysis. We introduce both discrete-time and continuous-time dynamic systems.

    A \textit{discrete-time dynamic system} is represented in the form $\mathbf{x}_{t+1} = \phi(\mathbf{x}_{t})$, where the state $\mathbf{x}_{t}\in M$ is defined at some time $t\in\mathbb{Z}$ on a compact $n$-dimensional manifold $M$ where $\phi: M \rightarrow M$ is a diffeomorphism, also called the iteration map, describing the evolution of the system with time. A \textit{continuous-time dynamic system} is a set of first-order ordinary differential equations which is represented in the form $\dot{\mathbf{x}} = v(\mathbf{x})$, where the state $\mathbf{x}\in M$ is defined at all time $t\in\mathbb{R}$, $M$ is a smooth compact $n$-dimensional manifold and $v$ is smooth ($\mathcal{C}^{2}$) vector field on $M$. The unique solution $\psi(\mathbf{x}, t): M \times \mathbb{R} \rightarrow M$ is the \textit{flow} generated by the vector field $v$, denoted as $\psi_{t}(\mathbf{x}) = \psi(\mathbf{x}, t)$ satisfying $\psi_{t_{1}}\circ\psi_{t_{2}} = \psi_{t_{1} + t_{2}}$ for $t \in \mathbb{R}$. Give a certain initial condition $\mathbf{x}_{0}\in M$, the solution to the dynamic system is the \textit{trajectory}, which is denoted as $\{\psi(\mathbf{x}_{t})\}$, for $t = 0,1,2,...$ or a curve $\{\psi_{t}(\mathbf{x}_{0})\}_{t\geq 0}$. Moreover, if the vector field $v$ or the iteration map $\phi$ is time-independent, the dynamic system is an \textit{autonomous dynamic system}. In this paper, our analysis is confined to autonomous systems. This focus is justified because a non-autonomous dynamical system can be transformed into an autonomous system by introducing the time variable as an additional dimension in the phase space. For example, the non-autonomous system $\dot{\mathbf{x}} = v(\mathbf{x}, t)$ is equivalent to the autonomous system in one higher dimension $\dot{\mathbf{x}} = v(\mathbf{x}, s)$, with $\dot{s} = 1$ where $s = t$.

    It is important to note that although the system is inherently continuous, practical constraints necessitate that measurements are taken at a regular sampling interval $T$. Given this sampling rate, the discrete flow denoted as $\xi_{T}: M \rightarrow M$, can be defined and characterized by the equations $\xi_{T}(\mathbf{x}_{t}) = \mathbf{x}_{t + T}$ and $\xi_{T}^{-1}(\mathbf{x}_{t})) = \mathbf{x}_{t-T}$ which associate each point state $\mathbf{x}\in M$ the vector $\xi_{T}(\mathbf{x})$. The sampling rate $T$ represents an integer multiple of the iteration step for discrete dynamic systems.

    In the theoretical framework of CCM, signals are causally linked if they share the common attractor manifold $\mathbb{A}$. The \textit{attractor} $\mathbb{A}$ is defined as a closed subspace on a smooth compact manifold $M$ with $dim(\mathbb{A})\leq n$ which satisfies the following three axioms \cite{strogatz2018nonlinear}:
    \begin{enumerate}
        \item $\mathbb{A}$ attracts an open set of initial conditions: There is an open set \(U\) containing $\mathbb{A}$ such that if the initial state $\mathbf{x}_{0} \in U$, then \[ \lim_{t\to\infty} dist(\mathbf{x}_{t}, \mathbb{A})=0 \] 
        \item Invariant set: If $\mathbf{x}_{t_{0}} \in \mathbb{A}$, then $\mathbf{x}_{t} \in \mathbb{A}$ for all $t\ge t_{0}$ 
        \item Minimal: No proper subset of $\mathbb{A}$ also satisfies these conditions.
    \end{enumerate} 
    The first two axioms ensure that points within the \textquotedbl basin of attraction\textquotedbl 
    are drawn towards $\mathbb{A}$, while the third establishes that every part of $\mathbb{A}$ is crucial. The existence of attractors has been studied for several decades, and it is established that every smooth, compact dynamical system possesses at least one attractor \cite{milnor1985concept}. Without loss of generality, we assume that the initial state $\mathbf{x}_{0}$ is contained within the attractor, ensuring that the trajectory remains within the attractor $\mathbb{A}$ during the periods of interest, such that $\mathbf{x}_{t}\in \mathbb{A} \subset \mathbb{R}^{N}$ for $t\geq 0$.

    Direct observation of the full state of the dynamic system is often infeasible. Instead, we have access to observations via a real-valued measurement function $h:\mathbb{A} \rightarrow \mathbb{R}$, producing the signals $\{s_{i}\}_{i\in\mathbb{N}} = \{h(\mathbf{x}_{i\cdot T})\}_{i}$, for $i = 0,1,2...$, where $\mathbf{x}_{0}$ is the initial state. A \textit{homeomorphism} between two manifolds $M_{1}$ and $M_{2}$ is a continuous bijection $f: M_{1} \rightarrow M_{2}$, where its inverse function $f^{-1}: M_{2} \rightarrow M_{1}$ is also continuous. Moreover, if the homeomorphism and its inverse are smooth, it is a diffeomorphism. An \textit{embedding} is a diffeomorphism from a manifold $M_{1}$ into another manifold $M_{2}$, defined as $f: M_{1} \rightarrow f(M_{1}) \subset{M_{2}}$. An important point is that embeddings are always injective and without self-intersections. Moreover, our goal is to find an embedding to reconstruct the attractor $\mathbb{A}$ from the signal $\{s_{i}\}_{i}$. Given the certain measurement function $h$, the following theorem forms the theoretical foundation for attractor reconstruction.  
    
    \begin{theorem}[Takens] \cite{takens2006detecting}
        Let \(M\) be an $n$-dimensional smooth manifold. If $v$ is a vector field on \(M\) with flow $\psi_{t}$ and $h$ is a measurement function on \(M\), then for generic choices of $v$ and $h$, the differential mapping $F_{h,m}: M \rightarrow \mathbb{R}^{m}$ of the continuous dynamic system into $\mathbb{R}^{m}$ is given by:
        \begin{equation}\label{derivative coordinate map}
            \begin{aligned}
            F_{h,m}(\mathbf{x}) = (h(\mathbf{x}), \frac{d}{dt}\Big|_{0}h(\psi_{t}(\mathbf{x})),...,\frac{d^{m-1}}{dt^{m-1}}\Big|_{0}h(\psi_{t}(\mathbf{x})))
            \end{aligned}
        \end{equation} 
        which is an embedding when $m = 2n + 1$, where $m$ is the embedding dimension, $\frac{d}{dt}\Big|_{0}$ means the derivatives are evaluated at $t= 0$ and the flow $\psi$ satisfies
        \begin{equation}
            \begin{aligned}
                \frac{d}{dt}\Big|_{0}\psi_{t}(\mathbf{x}) = v(\psi_{0}(\mathbf{x}))
            \end{aligned}
        \end{equation} 
        for every time $t \in \mathbb{R}$. 
    \end{theorem}
    \noindent The above theorem also holds for discrete dynamic systems with a diffeomorphism $\psi$ on a compact $n$-dimensional manifold $M$ and a measurement function $h$, for which the embedding is defined as Equation (\ref{delay-coordinate map}), where the value of the lag value $\tau$ is an integer multiple of the iteration size. The generic in this theory means that the differential mapping $F_{h,m}$ is an open and dense embedding in the set of all mappings under the measurement function $h$ and the flow $\psi_{t}$. The best way to understand this is regarding this theorem as a generalization of the Weak Whitney Embedding Theorem \cite{whitney1944self}. 
    \begin{theorem}[Weak Whitney Embedding]
        Every $n$-dimensional manifold M embeds in $\mathbb{R}^{2n+1}$. 
    \end{theorem}
    This theorem states that any manifold $M$ can be embedded in $\mathbb{R}^{2n+1}$ without self-intersections given an arbitrary mapping. Whitney proves that the optimal linear bound for the minimum embedding dimension is $2n$. Takens theorem demonstrates that the differential mapping (\ref{derivative coordinate map}) satisfies this condition, embedding the compact $n$-dimensional manifold $M$ into the reconstructed space $\mathbb{R}^{2n+1}$, even when considering finite discrete samples.

    For practical use, discrete versions of the differential mapping (\ref{derivative coordinate map}) are required when working with signals $\{s_{i}\}_{i}$ generated by the discrete flow $\xi_{T}$ with a specific sampling interval $T$. The most common approach is the delay-coordinate mapping $F_{h,\tau,m}(\mathbf{x}_{i \cdot \tau}): M\rightarrow\mathbb{R}^{m}$, which is defined as: 
    \begin{equation}\label{delay-coordinate map}
    F_{h,\tau,m}(\mathbf{x}_{i \cdot \tau}) = 
    \begin{bmatrix}
        h(\mathbf{x}_{i \cdot \tau}) \\
        h(\mathbf{x}_{(i-1) \cdot \tau}) \\
        \vdots \\
        h(\mathbf{x}_{(i-m+1) \cdot \tau})
    \end{bmatrix}
    =
    \begin{bmatrix}
        h(\mathbf{x}_{i \cdot \tau}) \\
        h(\xi_\tau^{-1}(\mathbf{x}_{i \cdot \tau})) \\
        \vdots \\
        h(\xi_\tau^{-m+1}(\mathbf{x}_{i \cdot \tau}))
    \end{bmatrix}    
    \end{equation}
    where the parameter $\tau = k\cdot T$, for $k\in\mathbb{Z}$ is the lag value, and $m$ is the embedding dimension. Theoretically, for minimal time delay $\tau$, a linear combination of coordinates can approximate the derivative such that the delay-coordinate mapping plays the same role as the differential mapping. The well-defined differential mapping is suitable for analytical purposes. This paper explores the properties of shadow manifolds reconstructed through differential mappings while implementing experiments using delay-coordinate mappings.

    In practical scenarios, the sampling rate $T$ of signals often cannot be small enough to accurately approximate the differential and higher-order differentials at the given point. However, by selecting an appropriate lag value $\tau$, the delay-coordinate mapping method can obtain the same result in reconstructing the shadow manifold using a sufficiently small $\tau$. Since chaotic dynamic systems consist of highly nonlinear and coupled variables, the signal obtained from the projection function, which serves as the measurement function, has the potential to recover information from other dimensions. The differential mapping method works by separating coupled information and projecting the observed data—via differentiation—in the direction of maximum linear independence, thereby isolating information about variables that are not directly observed. The critical part lies in accurately recovering information from the unknown dimensions using the observed data.

\begin{figure}[!ht]
        \centering
        \includegraphics[width=14cm]{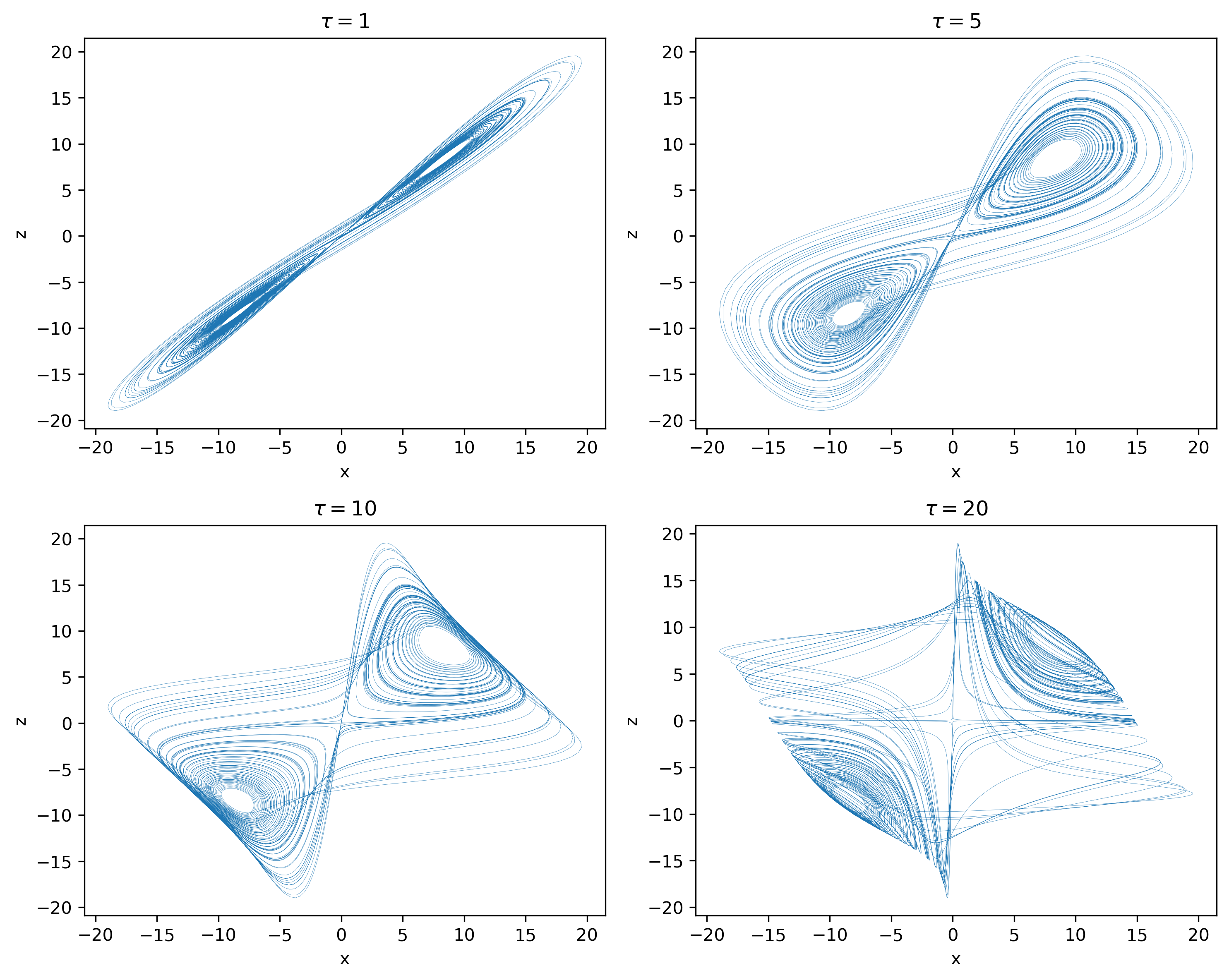}
        \caption{Sequential figures representing the changes in reconstructed shadow manifold $\mathcal{M}_{x}$ of Lorenz63 system with increasing lag value $\tau$.}
        \label{Increasing Lag Value}
    \end{figure}

    Thus, selecting the lag value $\tau$ plays a critical role in reconstructing the shadow manifold. If the lag value is suitable, the delay-coordinate mapping $\mathbf{F}_{h,\tau,n}(\mathbf{x}(t))$ is equivalent to the differential mapping $\mathbf{F}_{h,n}(\mathbf{x})$ under an affine transformation \cite{tsankov2004embeddings}, and plays as a diffeomorphism between the shadow manifold and the original attractor. However, the selection of the lag value is not only restricted by external factors like the sampling rate $T$ but also its intrinsic properties. For a continuous-time dynamic system with discrete flow, if $\tau$ is too small, the resulting vectors may be highly linear dependent and redundant, leading to a "squeezed" shadow manifold. Conversely, if $\tau$ is excessively large, the new coordinates may become essentially unrelated \cite{nichols2001attractor}, causing the shadow manifold to collapse. Based on the above analysis, we can observe that as $\tau$ increases, the shadow manifold undergoes a "stretch-and-fold" process, as depicted in Fig. \ref{Increasing Lag Value}. For convenience, we omit the sampling interval $T$ for $\tau$ such that the number of $\tau$ shown in this paper refers to the $k$ in the definition, indicating the number of times the sampling interval $T$, for example, $\tau = 5$ means $\tau = 5T$, where T is the sampling rate or the iteration steps for the discrete dynamic system.           

    Although the selection of lag value for delay-coordinate mapping is an open problem, several works have been done in this field \cite{tan2023selecting, martin2024robust}. The most widely-used method to choose the suitable lag value $\tau$ is the mutual information method \cite{kim1999nonlinear}. The basic idea is to calculate the mutual information between the system's observed values at different lag values and the original observed data and then select the first lag value at which the mutual information transitions from decreasing to increasing as the optimal $\tau$ since this represents the lag value that contains sufficient new information while still maintaining some correlation with the original data. This information-based method is theoretically intuitive, but the resulting values often do not correspond to the points at which the shadow manifold is fully stretched before collapsing. Furthermore, in cases where mutual information monotonically decreases with increasing lag value, this method does not work.

    According to Takens' theorem, a dynamic system can be embedded in an Euclidean space without self-intersection through any mapping. However, it should be noted that the ideal dimension for embedding an $n$-dimensional dynamic system is $n$. In many cases, the optimal dimension for the shadow manifold generated by differential embedding is also $n$. However, embedding an $n$-dimensional dynamic system into a higher dimension (greater than $n$) does not significantly affect the results of CCM, which is consistent with the results of tests using false nearest neighbors \cite{rhodes1997false}, as the redundant information introduced by the extra dimensions does not provide extra information about the original dynamic system and, therefore, does not impact the prediction. Therefore, when the dimension of the original dynamic system is known, the differential embedding method can be directly used to obtain a shadow manifold of the same dimension. In cases where information about the dimensions of the original dynamic system is lacking, false nearest neighbors are a practical approach for dimension selection.

    \begin{table}[H]
        \centering
        \begin{tabular}{ll}
            \hline
            \textbf{Notation} & \textbf{Meaning} \\
            \hline
            $\mathbf{x}_{t}$ & The state of dynamic system \\
            $M$ & Compact $n$-dimensional smooth manifold \\
            $\mathcal{M}$ & Reconstructed shadow manifold \\
            $h$ & Measurement function \\
            $\phi$ & Iteration map for discrete-time dynamic system \\
            $\mathbf{x}_t$ & State at time $t$, where $\mathbf{x}_t \in M$ \\
            $\psi_{t}$ & Flow generated by the vector field $v$ \\
            $v$ & Smooth vector field on $M$ \\
            $\psi(\mathbf{x}_t)$ & Trajectory of the system \\
            $\xi_T$ & Discrete flow with sampling rate $T$ \\
            $\mathbb{A}$ & Attractor \\
            $F_{h,m}$ & The differential mapping \\
            $\tau$ & Lag value for the delay-coordinate mapping $F_{h,\tau,m}$\\
            $F_{h,\tau,m}$ & The delay-coordinate mapping\\
            \hline
        \end{tabular}
        \caption{A list of mathematical notations in Section \ref{SectionTwo}.}
        \label{tab:notation}
    \end{table}

    \subsection{Convergent Cross Mapping}
    Here, we demonstrate how Takens' embedding theorem can be applied to detect causality. Given two-time series variables, \(X\) and \(Y\), which originate from two deterministic dynamic systems by using projection functions. Here we say \(X\) causes \(Y\) means the evolution of $\mathbf{y}_{t}$ depends on $\mathbf{x}_{t}$. The solid mathematical framework of CCM based on Takens' theorem is developed by Cummins et al. \cite{cummins2015efficacy} by interpreting unidirectional causality and bidirectional causality as projection and homeomorphism between two shadow manifolds separately. More precisely, if \(X\) unidirectionally causes \(Y\) (meaning \(Y\) does not cause \(X\)), then the joint dynamic system with an attractor $\mathbb{A}_{xy}$ can be constructed by concatenating \(X\) and \(Y\) as follows:   
    \begin{equation}
        \begin{aligned}
        \begin{bmatrix}
            x_{t+1} \\
            y_{t+1} 
            \end{bmatrix}
            =
            \begin{bmatrix}
            \psi_x(x_t) \\
            \psi_y(x_t, y_t)
            \end{bmatrix}
            = \psi_{xy}(x_t, y_t).
        \end{aligned}
    \end{equation} 
    The attractor $\mathbb{A}_{x}$ is then only a subset of the joint attractor $\mathbb{A}_{xy}$, which is identical to $\mathbb{A}_{y}$ and there is a noninjective projection $\pi_{yx}: \mathbb{A}_{y} \rightarrow \mathbb{A}_{x}$ between them. Consequently, there is also a noninjective map $\Pi_{yx}: \mathcal{M}_{y} \rightarrow \mathcal{M}_{x}$ from the shadow manifold $\mathcal{M}_{y}$ to $\mathcal{M}_{x}$ and the induced noninjective map $\Tilde{\Pi}_{y,x}$. In other words, we can reconstruct the shadow manifold $\mathcal{M}_{x}$ from $\mathcal{M}_{y}$, while the inverse is not guaranteed. Furthermore, if \(X\) causes \(Y\) and \(Y\) also causes \(X\), then the attractors $\mathbb{A}_{x}$ and $\mathbb{A}_{y}$ are equal to the joint attractor $\mathbb{A}_{xy}$. In this case, the shadow manifold $\mathcal{M}_{x}$ and $\mathcal{M}_{y}$ are diffeomorphic to the joint attractor. There exists a homeomorphism $f: \mathcal{M}_{x} \rightarrow \mathcal{M}_{y}$ between two shadow manifolds and an induced homeomorphism $f'$ between $\mathbb{A}_{x}$ and $\mathbb{A}_{xy}$. Consequently, each shadow manifold can be reconstructed from the other. The corresponding relationship of mapping is shown in Fig. \ref{CCM Mapping}.  
    \begin{figure}[h]
        \centering
        \resizebox{0.9\textwidth}{!}{ 
            \begin{minipage}{0.47\textwidth}
                \centering
                \begin{tikzpicture}[>=latex, node distance=2cm, thick]
                    \node (M1) at (0, 0) {$\mathcal{M}_{x}$};
                    \node (M2) at (4, 0) {$\mathcal{M}_{y}$};
                    \node (MU) at (2, 3) {$\mathbb{A}_{xy}$};
                    \draw[->] (MU) -- (M1) node[midway, left] {$F_{x,m}$};
                    \draw[->] (MU) -- (M2) node[midway, right] {$F_{y,m}$};
                    \draw[->, dashed] (M1) -- (M2) node[midway, below] {$f$};
                \end{tikzpicture}
                \caption*{(a)}
            \end{minipage}
            \hspace{0.02\textwidth}
            \begin{minipage}{0.47\textwidth}
                \centering
                \begin{tikzpicture}[>=latex, node distance=2cm, thick]
                    \node (M1) at (0, 0) {$\mathcal{M}_{x}$};
                    \node (M2) at (4, 0) {$\mathcal{M}_{y}$};
                    \node (MV) at (4, 3) {$\mathbb{A}_{xy}$};
                    \node (MU) at (0, 3) {$\mathbb{A}_{x}$};
    
                    \draw[->] (MU) -- (M1) node[midway, left] {$F_{x,m}$};
                    \draw[->] (MV) -- (M2) node[midway, right] {$F_{y,m}$};
                    \draw[->, dashed] (MU) -- (MV) node[midway, above] {$f'$};
                    \draw[->] (M1) -- (M2) node[midway, below] {$f$};
                \end{tikzpicture}
                \caption*{(b)} 
            \end{minipage}
        } 
        \resizebox{0.9\textwidth}{!}{ 
            \begin{minipage}{0.47\textwidth}
                \centering
                \begin{tikzpicture}[>=latex, node distance=2cm, thick]
                    \node (M1) at (0, 0) {$\mathcal{M}_{x}$};
                    \node (M2) at (4, 0) {$\mathbb{A}_{x}$};
                    \node (MV) at (4, 3) {$\mathbb{A}_{xy}$};
                    \node (MU) at (0, 3) {$\mathcal{M}_{y}$};
                    \draw[->] (MV) -- (MU) node[midway, above] {$F_{y,m}$};
                    \draw[->] (MV) -- (M2) node[midway, right] {$\pi_{yx}$};
                    \draw[->, dashed] (MU) -- (M1) node[midway, left] {$\Pi_{y,x}$};
                    \draw[->] (M2) -- (M1) node[midway, below] {$F_{x,m}$};
                \end{tikzpicture}
                \caption*{(c)} 
            \end{minipage}
            \hspace{0.02\textwidth}
            \begin{minipage}{0.47\textwidth}
                \centering
                \begin{tikzpicture}[>=latex, node distance=2cm, thick]
                    \node (M1) at (0, 0) {$\mathcal{M}_{x}$};
                    \node (M2) at (4, 0) {$\mathbb{A}_{x}$};
                    \node (MV) at (4, 3) {$\mathbb{A}_{xy}$};
                    \node (MU) at (0, 3) {$\mathcal{M}_{y}$};
                    \draw[->] (MV) -- (MU) node[midway, above] {$F_{y,m}$};
                    \draw[->, dashed] (MV) -- (M2) node[midway, right] {$\Tilde{\Pi}_{y,x}$};
                    \draw[->] (M2) -- (M1) node[midway, below] {$F_{x,m}$};
                    \draw[->] (MU) -- (M1) node[midway, left] {$\Pi_{y,x}$}; 
                \end{tikzpicture}
                \caption*{(d)} 
            \end{minipage}
        }
        \caption{First line: Bidirectional causation \(X\)$\Leftrightarrow$\(Y\). (a). Induced homeomorphism $f$ between $\mathcal{M}_{x}$ and $\mathcal{M}_{y}$. (b). Induced homeomorphism $f'$ between $\mathbb{A}_{x}$ and $\mathbb{A}_{xy}$. Last line: Unidirectional causation \(X\)$\Rightarrow$\(Y\) when $f$ exists. (c). Induced noninjective projection $\Pi_{y,x}$ between $\mathcal{M}_{y}$ and $\mathcal{M}_{x}$. (d). Induced noninjective projection $\Tilde{\Pi}_{y,x}$ when $\Pi_{y,x}$ exists. }\label{CCM Mapping}
    \end{figure}
    
    Based on the above illustration, Sugihara et al. proposed the CCM algorithm to infer causal links between two signals. It aims to do this by using $k$-nearest-neighbor regression on one shadow manifold to make predictions of another and evaluating the forecasting skill via the Pearson correlation coefficient. Given two time series $\{x_{1}, x_{2},..., x_{l}\}$ and $\{y_{1}, y_{2},..., y_{l}\}$ from two time series variables \(X\) and \(Y\) with finite length $l$, CCM tests the causality from $X \Rightarrow Y$ in the following three steps: 
    
    \begin{enumerate}
        \item Reconstruct two shadow manifolds $\mathcal{M}_{x}$ and $\mathcal{M}_{y}$ using the delay-coordinate mappings $F_{h,\tau,n}$.
        \item For each point $\mathbf{y}_{t}$ in the shadow manifold $\mathcal{M}_{y}$, collect its $n+1$ nearest points $\mathbf{y}_{t_{i}}, i=1,2,...,n+1$. Estimate $\hat{{x}}_{t}|\mathcal{M}_y$ using locally weighted averages of the corresponding ${x}_{t_{i}}$ values from the original time series \(X\), calculated as:
        \begin{align}
            \hat{{x}}_{t}|\mathcal{M}_{y} = \sum_{i=1}^{n+1} w_{i}{x}_{t_{i}}
        \end{align}
        where the weights $w_{i}$ depends on the distance between $\mathbf{y}_{t}$ and its $i$-th closest neighbor in $\mathcal{M}_{y}$, computed as:
        \begin{align}
            w_{i} = \frac{u_{i}}{\sum_{j=1}^{n+1} u_{j}}
        \end{align} 
        with $u_{i}$:
        \begin{align}
            u_{i} = \rm{exp}\{\frac{-d[\mathbf{y}_{t},\mathbf{y}_{t_{i}}]}{d[\mathbf{y}_{t},\mathbf{y}_{t_{1}}]}\}
        \end{align}
        where $d[\mathbf{y}_{t},\mathbf{y}_{t_{1}}]$ represents the Euclidean distance between two vectors $\mathbf{y}_{t}$ and $\mathbf{y}_{t_{1}}$. Similarly, cross mapping in the other direction is defined by swapping variables.  
        \item Calculate the Pearson correlation coefficient between $\hat{{x}}_{t}|\mathcal{M}_y$ and ${x}_{t}$, and its absolute value is the corresponding forecast skill $\rho_{xy}$.  
        Repeat 1 and 2 sequentially as the library length $l$ increases, and observe whether the value of $\rho_{xy}$ converges to a nonzero value.    
    \end{enumerate}
    
    If \(X\) causes \(Y\), then as the library length $l$ increases, the shadow manifold becomes denser, allowing the estimator $\hat{{x}}_{t}|\mathcal{M}_{y}$ converge to ${x}_{t}$. Consequently, the Pearson correlation between the ground truth and the predictor gradually converges. In other words, as we increase the length of the signal if the predictability of newly added points aligns with the previous points, it suggests the presence of causality rather than a mere statistical coincidence. The Pearson correlation $\rho$ value may converge to a high plateau; otherwise, the cause-and-effect relationship does not exist.

    Moreover, it is important to note that in cases of strong coupling, for instance, when \(X\) significantly influences \(Y\), it may exhibit behaviors that are not independent of X. In such scenarios, the causal links may appear reversible, leading to an erroneous interpretation of unidirectional causation as bidirectional. This phenomenon is known as the problem of generalized synchrony. Further research on detecting generalized synchrony within the CCM framework is discussed in \cite{ye2015distinguishing}.

    \subsection{Related Work about CCM}

    Several factors, such as noisy data, irregular or sporadic sampling, and improper selection of embedding parameters, significantly influence the correctness of CCM. Here, we introduce some related work as refinements of CCM. A key limitation of CCM is its sensitivity to missing data, necessitating long, uninterrupted data for regular sampling \cite{de2020latent}. Latent-CCM has been developed to address this in repeated short, sporadic time series. It leverages a Neural ODE model to learn the inherent dynamics of the data, a method that provides more state space information than the multi-spatial CCM \cite{clark2015spatial}. The TD-CCM method, employing denoising techniques like wavelet methods and empirical mode decomposition, has been proposed in noisy data and trending noise signals.

    Furthermore, Feng \cite{feng2019detecting} introduced a Bayesian version of CCM, using a deep Gaussian process based on Cross Mapping Smoothness (CMS) \cite{ma2014detecting}. This approach shifts from the traditional online k-nearest-neighborhood approximation to a deep Gaussian process, underlining the correlation between map smoothness and the strength of time series relationships. It also includes a method for dimension reduction of the shadow manifold using the Gaussian Process Latent Model. The CMS aims to detect causal relations in short-term time series by equating map smoothness with causal detection, using a Radial Basis Function Network for training. While theoretically sound, this method struggles with differentiating whether the neural network's fitting ability or significant training errors contribute to the model's performance, making model selection crucial. Moreover, CMS's effectiveness for short-term data may be limited due to potential overfitting issues, as smoothness is a necessary but insufficient condition for cross mapping.

    Finally, Bulter \cite{butler2023causal} highlights the importance of auto-predictability and recurrence as prerequisites for employing CCM effectively. Autopredictability ensures the deterministic nature of the system for reliable prediction based on its history, while recurrence validates the presence of neighborhoods on the attractor manifold. Fulfilling these criteria is crucial for successful state space construction, and additional experiments \cite{cobey2016limits, krakovska2020implementation} are recommended to test the statistical significance of CCM results.

    \section{Problem Formulation}\label{SectionThree}
    
    For a complex system with interacting variables connected, there are two kinds of causality: direct and indirect. A prototypical example of an indirect causal link is the variable \(Y\) within the Lorenz63 system, which acts as an interacting variable, introducing indirect causal links between \(X\) and \(Z\). CCM can indiscriminately detect such indirect causality due to the causation transitivity \cite{leng2020partial, ghouse2021inferring}. Although CCM cannot distinguish those two types of links, \cite{leng2020partial} provides the Partial Cross Mapping method to eliminate indirect causal influence. Based on this, we do not discuss whether the causality detected by CCM is indirect or direct in this paper. This section shows that CCM may not provide the correct causal links for the Lorenz63 system. Notably, this counterexample is not a statistical anomaly. We have identified similar counterexamples originating from various systems; for example, the Chen $\&$ Ueta dynamic system \cite{chen1999yet} wherein CCM fails to detect cross maps introduced by the interacting variable \(Y\).

\begin{figure}[h]
    \centering
    \includegraphics[width=10cm]{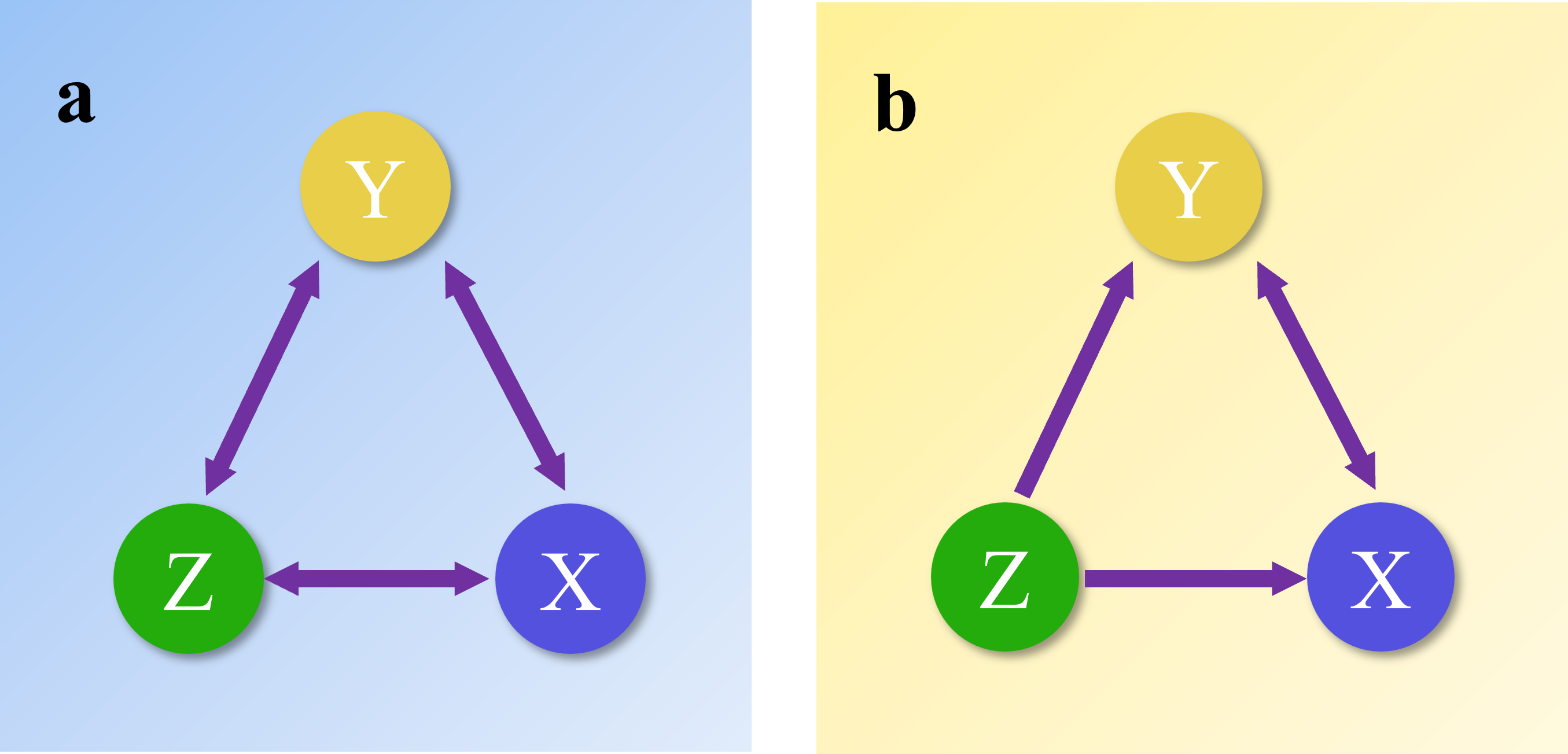}
    \caption{True causality versus CCM's result. a. Graphical model of causal relations within the Lorenz63 system. b. Graphical model of causal relations obtained by CCM. The Lorenz63 system is simulated using the fourth-order Runge-Kutta method, where the temporal domain is $t\in[0,100]$, the time step for discretization is $\Delta t=0.01$ which could be regarded as the sampling rate $T$ for continuous-time systems, the initial condition is $\mathbf{x}_{0} = [1,1,1]$, and $\tau = 9$, $n =3$ for CCM.}
    \label{CCM4Lorenz63}
    \end{figure}

\begin{figure}[h]
    \centering
    \includegraphics[width=16cm]{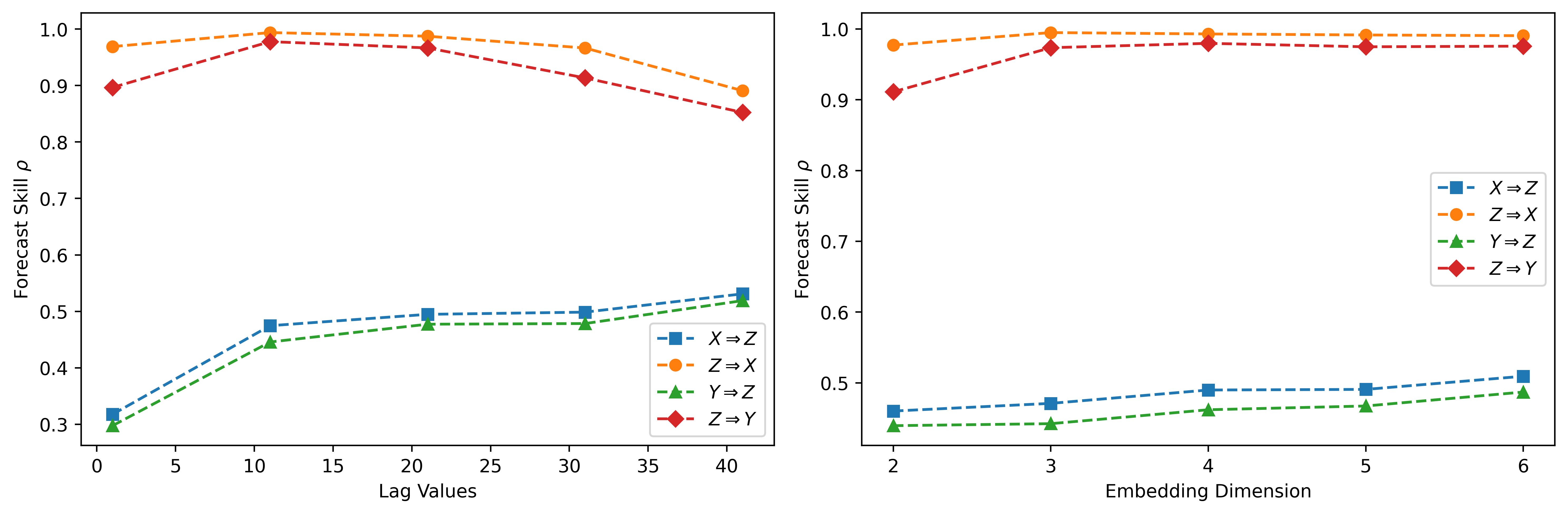}
    \caption{Influence of embedding parameters. Left. The influence of lag value $\tau$ when $n =3$. Right. The influence of embedding dimension $n$ when $\tau = 9$. }
    \label{Parameter Test}
    \end{figure}

    The classical chaotic Lorenz63 system is:
    \begin{equation}
        \begin{aligned}\label{LorenzFormulation}
        & \frac{dx}{dt} = \sigma(y-x)\\
        & \frac{dy}{dt} = \rho x - y - xz\\
        & \frac{dz}{dt} = xy - \beta z
        \end{aligned}
    \end{equation}
    where the parameters are defined as $\sigma = 10, \rho = 28$, and $\beta = \frac{8}{3}$ with the initial condition $\mathbf{x}_{0}=[1,1,1]$. In the Lorenz63 system, \(X\) influences both \(Y\) and \(Z\), \(Y\) affects both \(X\) and \(Z\), and \(Z\) affects \(Y\) and \(X\). However, as shown in Fig. \ref{CCM4Lorenz63}, CCM mistakes bidirectional causalities $Y \Leftrightarrow Z$ and $X \Leftrightarrow Z$ as unidirectional $Z \Rightarrow Y$ and $Z \Rightarrow X$. Multiple experiments have also been carried out to avoid the influence of parameters to test the influence of the value of $\tau$ and $n$ shown in Fig. \ref{Parameter Test}. 
    
    The fluctuation of lag values in a reasonable range may only fluctuate the convergent score of CCM in a small range, thus persisting in the same casual direction. However, the condition of the embedding dimension is different. The score is almost stable when the embedding dimension reaches the suitable number ($n=3$ in this case) because information obtained by increasing the dimension is redundant and never provides additional information for improving predicting skills. The failure in the Lorenz63 system poses a significant query: Is this a statistical anomaly? Further investigation should be done to find more similar counterexamples. Moreover, such counterexamples should be distinct from the Lorenz63 system, neither transformable into it nor derivable from it through linear or nonlinear coordinate transformations (i.e., diffeomorphisms).

    Besides the Lorenz63 system, we also discover numerous counterexamples for which CCM fails to detect the correct causal links. This section focuses primarily on three examples: the Chen \& Ueta system, the Burke \& Shaw system \cite{shaw1981strange}, and the three-scroll chaotic attractor \cite{li2008three}. Generally speaking, all three systems are Lorenz-like systems since they share some common properties with the Lorenz63 dynamic system \cite{letellier2023lorenz}, for example, similar butterfly attractors like the Lorenz63 attractor. However, they are topologically different and nonequivalent to the Lorenz63 system. Additional examples are elaborated in Section \ref{SectionFive}.

\begin{figure}[h]
        \centering
        \includegraphics[width=16cm]{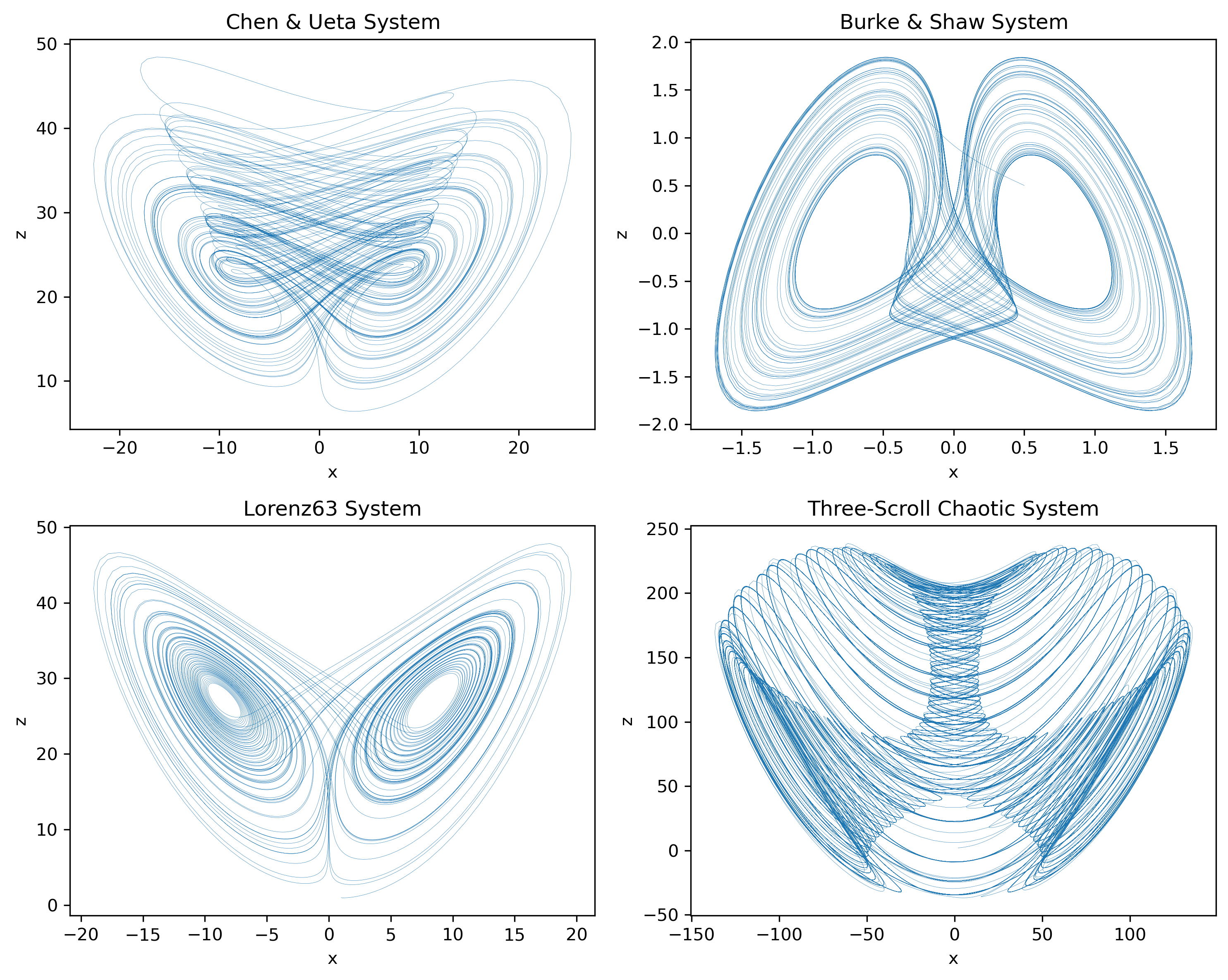}
        \caption{x-z plane projections of attractors. a. Chen $\&$ Ueta system, b. Lorenz63 system, c. Burke $\&$ Shaw system, d. three-scroll chaotic system.}
        \label{LorenzAndChen}
    \end{figure}

    The first example is the Chen $\&$ Ueta dynamic system; its mathematical representation is formulated as follows:
    \begin{equation}
    \begin{aligned}
        &\frac{dx}{dt} = \alpha(y-x)\\
        &\frac{dy}{dt} = (\gamma-\alpha)x - xz + \gamma y\\
        &\frac{dz}{dt} = xy - \beta z
    \end{aligned}
    \end{equation}
    where the initial condition is set as $x_{0} = [-10,0,37]$. This system, characterized by its chaotic nature and multiple attractors, is shown in Fig. \ref{LorenzAndChen}, with parameters $\alpha = 35$, $\beta = 3$, and $\gamma = 28$.  
    Despite its mathematical similarity to the Lorenz63 system, the Chen $\&$ Ueta system exhibits topological differences, primarily in the number of equilibria: three in this case, as opposed to two in the Lorenz63. The causal structure of this system is the same as the Lorenz63. Another Lorenz-like system, the Burke \& Shaw system, is described as follows:
    \begin{equation}\label{Burke and Shaw}
    \begin{aligned}
        &\frac{dx}{dt} = -\alpha(x+y)\\
        &\frac{dy}{dt} = -y - \alpha xz\\
        &\frac{dz}{dt} = \beta + \alpha xy
    \end{aligned}
    \end{equation}
    with parameters $\alpha = -10.0$, $\beta = 4.272$ and initial state $\mathbf{x}_{0} = [0.5, 0.5, 0.5]$. The attractor produced by equation (\ref{Burke and Shaw}) is also topologically different from the one of the Lorenz63 system \cite{letellier1996evolution}.

    We simulated both systems using the fourth-order Runge-Kutta method. For the Chen $\&$ Ueta system, the temporal domain is $t\in[0,50]$, and $\Delta t=0.005$. Here, we choose a smaller $\Delta{t}$ because this system has a larger vector field, and a small step size allows for better illustration and simulation. The lag value and embedding dimension we optimized of the delay coordinate map as $\tau = 20$ and $n = 3$. For the Burke \& Shaw system, the trajectories are simulated in the temporal domain $t \in [0,100]$ with $\Delta t = 0.01$, and the embedding parameters we optimized are $\tau = 10$ and $n = 3$.

    The last counterexample is the three-scroll chaotic system:
    \begin{equation}
    \begin{aligned}
        &\frac{dx}{dt} = 40(y-x)+0.16xz\\
        &\frac{dy}{dt} = 55x + 20y -xz\\
        &\frac{dz}{dt} = -0.65x^{2} +xy + \frac{11}{6}z
    \end{aligned}
    \end{equation}
    with the initial point $\mathbf{x}_{0} = [2, 2, 2]$, temporal domain $t = [0,150]$, $\Delta t = 0.0015$, and $\tau = 20$, $n = 3$. This dynamic system generates a toroidal attractor bounded by a genus-3 torus, which is shown in Fig. \ref{LorenzAndChen}. Six causal links exist within this chaotic system as $X \Leftrightarrow Y$, $X \Leftrightarrow Z$, and $Y \Leftrightarrow Z$. Causal results of these three systems obtained through CCM are listed in Table \ref{Details of Lorenz-like system}, where we can see that CCM fails to provide the correct result. These instances suggest that the failure observed in the Lorenz63 is not merely an anomaly but indicative of an underlying commonality that warrants further investigation.                       
    
    \section{Explanation and Proposed Solution}\label{SectionFour}
    
    In the previous section, we proposed that CCM may provide incorrect causal links for several cases. Here, we explain the reason and extend it to general cases. We point out that for the $n$-fold dynamic system, which is equivariant under a nontrivial cyclic group $\mathbb{C}_{n}$ of order $n$, the differential mapping of the invariant variable is a $n$-to-one mapping rather than an embedding. Consequently, when there is a bidirectional causation between invariant and non-invariant variables, there is no longer a homeomorphism between two shadow manifolds, misleading a wrong causality between two variables. To address this problem, we propose a practical method for the case when the dynamic system is equivariant under the order two cyclic group $\mathbb{C}_{2}$.

    \subsection{Preliminary}
    
    Here, we briefly introduce some basic definitions and properties that will be used later. \textit{Group} theory is the mathematical tool to study symmetry. For a finite group $G$, the \textit{order} of the group refers to the number of elements it contains, denoted by $|G|$. A \textit{group action} \texttt{"$\cdot$"} of a group $G$ on $\mathbb{R}^{n}$ is defined as a map $G \times \mathbb{R}^{n} \rightarrow \mathbb{R}^{n}$ which satisfies: (a). $e\cdot x = x$ for all $x \in \mathbb{R}^{n}$, where $e$ is the identity element of the group $G$ and (b). $(gh)\cdot x = g\cdot(h\cdot x)$ for $x \in \mathbb{R}^{n}$ and $g, h \in G$. Here $g\cdot x$ means the image of $x$ under the action of $g \in G$.

    A dynamic system $\dot{\mathbf{x}} = v(\mathbf{x})$ with flow $\psi_{t}(\mathbf{x})$ is \textit{$|G|$-fold symmetric} or \textit{G-equivariant} if 
    \begin{equation}\label{Symmetric Property}
    \begin{aligned}
           g\cdot\dot{\mathbf{x}} = g\cdot v({\mathbf{x}}) = v(g\cdot{\mathbf{x}}).
    \end{aligned}
    \end{equation}
    for every $g\in G$ and $\mathbf{x}\in\mathbb{R}^{n}$. For every $g \in G$, the \textit{fundamental domain} $\mathcal{D}_{g}\in\mathbb{R}^{n}$ under the group $G$ is a representative region such that the action of $G$ tessellates the entire space without overlaps. Formally, it satisfies: 
    \begin{equation}
    \begin{aligned}
            g_{1}\cdot \mathcal{D}_{g_{1}} \cap g_{2}\cdot \mathcal{D}_{g_{2}} \neq \emptyset \Rightarrow g_{1} = g_{2}
    \end{aligned}
    \end{equation}
    for $g_{1}, g_{2} \in G$. For example, the fundamental domain of the Lorenz63 system can be chosen as $\mathcal{D}_{R_{z}} = \{(x,y,z) | x\geq 0\}$. If $G$ is a reflection group, then the dynamic system is \textit{reflection symmetric}. If $G$ is the cyclic group $\mathbb{C}_{n}$, then the dynamic system is \textit{rotation symmetric}. If the dynamic system is equivariant under the inversion group $G = \{e, P\}$ where $P$ denotes the inversion operation as $P: (x_{1}, x_{2},..., x_{n}) \rightarrow (-x_{1}, -x_{2},...,-x_{n})$, then the dynamic system is \textit{inversion symmetric}.

\begin{figure}[h]
        \centering
        \includegraphics[width=8cm]{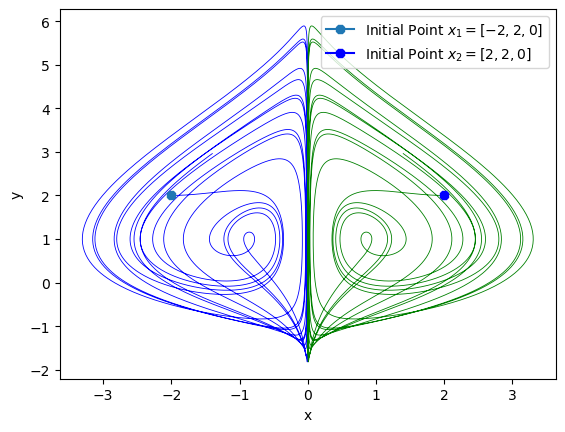}
        \caption{A symmetric pair of attractor "kissing" of the reflection equivariant system, when a = 0.7 and two initial values $\mathbf{x_{1}} = [-2, 2, 0]$ and $\mathbf{x_{2}} = [2, 2, 0]$. Two parts of trajectories can become arbitrarily close but never cross the plane $x = 0$.}
        \label{kissing}
    \end{figure}

    There are several crucial factors we should notice for a symmetric attractor. Firstly, the group $G$ should satisfy certain algebraic conditions when the reflection planes exist such that the $G$-equivariant dynamic system may generate a symmetric attractor \cite{ashwin1994symmetry, field1996symmetric}. Moreover, not all symmetric attractors are connected because of the existence of the reflection plane of the symmetry group. The system gives a typical example of this \cite{sprott2014simplest}:     
    \begin{equation}\label{simplest reflection}
    \begin{aligned}
            \dot{x} &= x - xy\\
            \dot{y} &= z\\
            \dot{z} &= -y-az+x^{2}
    \end{aligned}
    \end{equation}
    which is $G$-equivariant under the group $G = \{e, \sigma_{x}\}$, where $\sigma_{x}: (x,y,z) \rightarrow (-x,y,z)$, producing the so-called "kissing" attractor shown in Fig. \ref{kissing}.   
    It is evident that the symmetric pair of strange attractors nearly touch each other on opposite sides of the reflection plane $x=0$, as $\dot{x} = 0$ in this plane. Consequently, a trajectory originating from an initial point in one part of its disconnected symmetric attractor may not exhibit symmetry behavior. Based on the above analysis, when the symmetry group has a reflection plane, as in such cases, the corresponding trajectory may not show the expected symmetry. In this paper, we focus on investigating the impact of the symmetry of rotation symmetric systems with no reflection plane on detecting causality by CCM.

    The following are two useful results for differential mapping. For $g\in G$, if the function $u$ satisfying:
    \begin{equation}
    \begin{aligned}
           u(g\cdot\mathbf{x}) = \pm u(\mathbf{x}),
    \end{aligned}
    \end{equation}
    we say that $u$ defines the \textit{(even or odd) parity} under $g$. For a G-equivariant dynamic system $\dot{\mathbf{x}} = v(\mathbf{x})$ with flow $\psi_{t}(\mathbf{x})$, by the definition of the derivative, we have: 
    \begin{equation}
    \begin{aligned}
           \frac{d}{dt}\Big|_{0}u[\psi_{t}(\mathbf{x})] = \lim_{t=0}\frac{u[\psi_{t}(\mathbf{x})] - u(\mathbf{x})}{t} = \lim_{t=0}\frac{u[\mathbf{x} + tv(\mathbf{x})] - u(\mathbf{x})}{t} 
    \end{aligned}
    \end{equation} 
    where $\frac{d}{dt}\Big|_{0}$ denotes that the derivative is evaluated at time $t = 0$, and $\psi_{t}(\mathbf{x})$ represents the flow of the system. By using the property in (\ref{Symmetric Property}), and transforming $\mathbf{x}$ to $g\cdot\mathbf{x}$, we obtain:
    \begin{equation}\label{derivative}
    \begin{aligned}
           \frac{d}{dt}\Big|_{0}u[\psi_{t}(g\cdot\mathbf{x})] &= \lim_{t=0}\frac{u[\psi_{t}(g\cdot\mathbf{x})] - u(g\cdot\mathbf{x})}{t}\\
           &= \lim_{t=0}\frac{u[g\cdot\mathbf{x} + tv(g\cdot\mathbf{x})] - u(g\cdot\mathbf{x})}{t} \\
           &= \lim_{t=0}\frac{u[g\cdot(\mathbf{x} + tv(\mathbf{x}))] - u(g\cdot\mathbf{x})}{t} \\ 
           &= \frac{d}{dt}\Big|_{0}u[g\cdot\psi_{t}(\mathbf{x})] 
    \end{aligned}
    \end{equation} 
    Thus, it implies that the time derivative of $u$ also preserves the same parity. By mathematical induction, we can conclude that the $n$-th derivative of $u$ inheres the same parity.

    Through the differential mapping $F_{h,n}$ defined as (\ref{derivative coordinate map}), the shadow manifold is constructed from consecutive derivatives. The corresponding vector field of the shadow manifold is:
    \begin{equation}
    \begin{aligned}
           V^{i} =  \frac{d}{dt}\Big|_{0}F_{h,n}^{i}(\psi_{t}(\mathbf{x}))\\
    \end{aligned}
    \end{equation} 
    where $F_{h,n}^{i}$ and $V^{i}$ denote the $i$-th component of $F_{h,n}(\mathbf{x})$ and $V$ separately, and we have the following rules:
    \begin{equation}
    \begin{aligned}
           \frac{d}{dt}\Big|_{0} F_{h,n}^{2}(\psi_{t}(\mathbf{x})) = F_{h,n}^{3}(\mathbf{x}).\\
    \end{aligned}
    \end{equation} 
    Since $V^{1} = F^{2}$, by mathematical induction, we have the general rule:
    \begin{equation}\label{vector field}
    \begin{aligned}
           V^{i}(\mathbf{x}) = \frac{d}{dt}\Big|_{0}F_{h,n}^{i}(\psi_{t}(\mathbf{x})) = F^{i+1}_{h,n}, \quad i<n
    \end{aligned}
    \end{equation} 
    Thus, the $n$-dim shadow manifold always has the analytic canonical form:
    \begin{equation}
    \begin{aligned}
    V^{1} &= \frac{d}{dt}\Big|_{0} F_{h,n}^{1} = F_{h,n}^2, \\
    V^{2} &= \frac{d}{dt}\Big|_{0} F_{h,n}^{2} = F_{h,n}^3, \\
    &\vdots \\
    V^{n-1} &= \frac{d}{dt}\Big|_{0} F_{h,n}^{n-1} = F_{h,n}^n, \\
    V^{n} &= \frac{d}{dt}\Big|_{0} F_{h,n}^{n} = f(F_{h,n}^1, \dots, F_{h,n}^n).
    \end{aligned}
    \end{equation}
    for some function $f$, which is the only unknown term.

    The following theorem explains the preservation of symmetry in the reconstructed system when using differential mapping to reconstruct a symmetric system from a single observation function \cite{cross2010equivariant}.     
    \begin{theorem}[Cross]\label{Cross}
        A differential reconstruction of any nonlinear dynamic system preserves, at most, a two-fold symmetry. 
    \end{theorem}
    Specifically, there are only two possibilities for reconstructed shadow manifolds: either (a). The shadow manifold has no symmetry, and (b). The shadow manifold is equivariant under the inversion group $G = \{e, P\}$.

    \subsection{Explanation}
    
    We begin with a few key observations. First, all three dynamic systems are two-fold symmetric, specifically three-dimensional ones that rotate around the $z$-axis. In such systems, the differential mapping $F_{z,n}$ generates two symmetric shadow manifolds, $\mathcal{M}_{x}$ and $\mathcal{M}_{y}$, along with a nonsymmetric one, $\mathcal{M}_{z}$. Second, all trajectories of these shadow manifolds are connected. These observations lead us to investigate the properties of symmetric systems under differential mappings. Here, we give the following remark as an explanation for the general case.    

    \textbf{Proposition 1 .} \textit{If the dynamic system $\mathbf{\dot{x}} = v(\mathbf{x})\in \mathbb{R}^{n}$ is symmetric under the order-two cyclic group $\mathbb{C}_{2} = \{e, R_{x_{n}}(\pi)\}$, where $R_{x_{n}}(\pi): (x_{1}, x_{2}, ..., x_{n})\rightarrow (-x_{1},-x_{2},...,x_{n})$, CCM may mislead the bidirectional causality $x_{n}\Leftrightarrow x_{i}$ between the invariant variable $x_{n}$ and $x_{i}, i\neq n$ as unidirection relationship $x_{n} \Rightarrow x_{i}$.} 
    \renewcommand\qedsymbol{$\blacksquare$}
    
    \begin{proof}
    It is sufficient to show that there is a noninjection mapping $\Pi_{x_{i}x_{n}}: \mathcal{M}_{x_{i}} \rightarrow \mathcal{M}_{x_{n}}$ between $\mathcal{M}_{x_{i}}$ and $\mathcal{M}_{n}$. Since the dynamic system is two-fold symmetric under $\mathbb{C}_{2}$, then when the measurement function $h$ is taken to be the $x_{n}$-coordinate projection, it defines even parity as $h(R_{x_{n}}(\pi)\cdot\mathbf{x}) = h(\mathbf{x})$. Using the results (\ref{derivative}) and (\ref{vector field}), the differential mapping $F_{x_{i},n}$ and the new vector field also inherit the same even parity as $h$. Therefore, the $F_{x_{n},n}$ is noninjective and is two-to-one at the nicest case. Because of the bidirectional causality $x_{n} \Leftrightarrow x_{i}$, we can find the induced homeomorphism $f: \mathbb{A}_{x_{n}}\rightarrow\mathbb{A}_{x_{i}x_{n}}$ and $f': \mathbb{A}_{x_{i}x_{n}}\rightarrow\mathbb{A}_{x_{n}}$. Then the homeomorphism $\Tilde{f}: \mathcal{M}_{x_{i}}\rightarrow \mathbb{A}_{x_{n}}$ can be obtained by $\Tilde{f}: f'\circ F_{x_{i},n}^{-1}$, since $F_{x_{i},n}$ is a diffeomorphism. So the map  
    \begin{equation}
    \begin{aligned}
        \Pi_{x_{i}x_{n}} :&= F_{x_{n},n} \circ \Tilde{f}\\
                          &= F_{x_{n},n} \circ f'\circ F_{x_{i},n}^{-1}
    \end{aligned}
    \end{equation}
    constructed from Fig. \ref{Proof} is the desired noninjective mapping since $F_{x_{n},n}$ is noninjective. 
    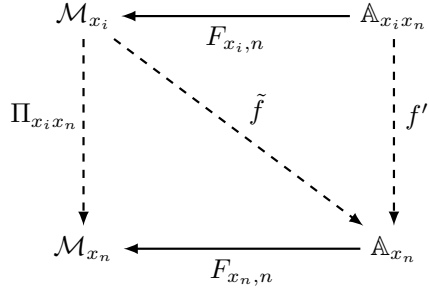
\begin{figure}[h]
        \centering
        \resizebox{0.5\textwidth}{!}{ 
            \begin{minipage}{0.47\textwidth}
                \centering
                \begin{tikzpicture}[>=latex, node distance=2cm, thick]
                    \node (M1) at (0, 0) {$\mathcal{M}_{x_{n}}$};
                    \node (M2) at (4, 0) {$\mathbb{A}_{x_{n}}$};
                    \node (MV) at (4, 3) {$\mathbb{A}_{x_{i}x_{n}}$};
                    \node (MU) at (0, 3) {$\mathcal{M}_{x_{i}}$};
                    \draw[->, dashed] (MU) -- (M2) node[midway, right, yshift=0.3cm] {$\Tilde{f}$};
                    \draw[->] (M2) -- (M1) node[midway, below] {$F_{x_{n},n}$};
                    \draw[->] (MV) -- (MU) node[midway, below] {$F_{x_{i},n}$};
                    \draw[->, dashed] (MV) -- (M2) node[midway, below, xshift=0.3cm, yshift=0.5cm] {$f'$};
                    \draw[->, dashed] (MU) -- (M1) node[midway, below, xshift=-0.5cm, yshift=0.5cm] {$\Pi_{x_{i}x_{n}}$};
                \end{tikzpicture}
            \end{minipage}
        }
        \caption{Induced projection $\Pi_{x_{i}x_{n}}: \mathcal{M}_{x_{i}}\rightarrow \mathcal{M}_{x_{n}}$, homeomorphism $f': \mathbb{A}_{x_{i}x_{n}}\rightarrow \mathbb{A}_{x_{n}}$ and $\Tilde{f}: \mathcal{M}_{x_{i}}\rightarrow \mathbb{A}_{x_{n}}$. }\label{Proof}
    \end{figure}
    \end{proof}


    Next, we use the Lorenz63 and the Burke $\&$ Shaw systems to show concrete examples. Both systems are two-fold symmetric system under the cyclic group $\mathbb{C}_{2} = \{e, R_{z}(\pi)\}$. For the Lorenz63 system, when the measurement function is taken to be the $x$-coordinate projection which defines odd parity under the action of $R_{z}(\pi)$, the corresponding differential mapping $F_{x, 3}: (x,y,z) \rightarrow (u,v,w)$ is as follows: 
    \begin{equation}
    \begin{aligned}
             u &= x \\
             v &= \sigma(y-x)  \\
             w &= \sigma[(\rho+\sigma)x - (\sigma + 1)y - xz] \\
    \end{aligned}
    \end{equation}
    where we use $(u, v, w)$ to denote new coordinate to avoid confusion. The new vector field of the $x$-induced shadow manifold also referred to as the induced Lorenz system, is:  
    \begin{equation}\label{induced Lorenz}
    \begin{aligned}
             \dot{u} &= v \\
             \dot{v} &= w  \\
             \dot{w} &= \beta\sigma(\rho-1)u - \beta(\sigma+1)v - (1+\beta+\sigma)w - u^{2}v - \sigma u^{3} + \frac{v}{u}(w+(1+\sigma)v) \\
    \end{aligned}
    \end{equation}
    where the domain is $\{(u,v,w)\in\mathbb{R}^{3}| u\neq 0\}$.

\begin{figure}[ht]
        \centering
        \includegraphics[width=14cm]{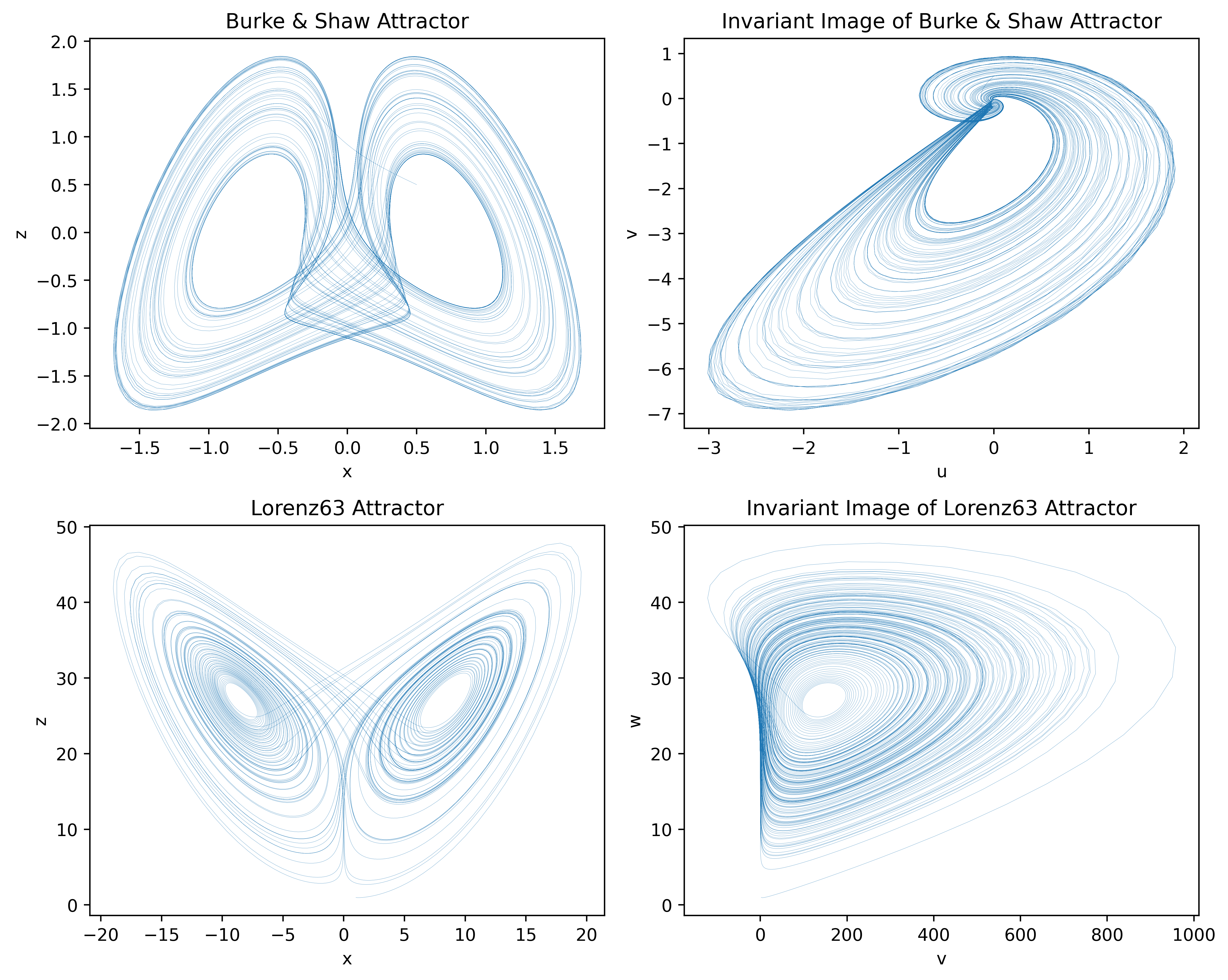}
        \caption{Left: Strange attractors generated by the Lorenz63 system (\ref{LorenzFormulation}) and Burke \& Shaw system (\ref{Invariant Burke and Shaw}). Right: Invariant images of the Lorenz63 and Burke \& Shaw systems.}
        \label{Invariant Image}
    \end{figure}

    The same process applies to obtain the invariant system generated by $F_{z, 3}$ for the Burke \& Shaw equation (\ref{Burke and Shaw}) as: 
    \begin{equation}\label{Invariant Burke and Shaw}   
    \begin{aligned}
             \dot{u} &=  -(\alpha + 1)u - \alpha(1 - w)v + (1 - \alpha)\rho \\
             \dot{v} &=  \alpha(1-w)u - \alpha(1+w)\rho - (\alpha + 1)v \\
             \dot{w} &=  \frac{\alpha}{2}v + \beta 
    \end{aligned}
    \end{equation}
    where $\rho = \sqrt{u^{2}+v^{2}}$. The invariant image of Lorenz63 and Burke \& Shaw is shown in Fig \ref{Invariant Image}.

    We have systematically examined the structure of two-fold rotation symmetric systems in $\mathbb{R}^{3}$, the reconstruction of their shadow manifolds, and the characteristics of the differential mapping $F_{h,3}$. This analysis explains why CCM fails to provide a correct causal graph consistent with other counterexamples. Specifically, the variable that defines even parity and remains invariant under the representation of the order-two group lacks sufficient information about the attractor's symmetry, making it unable to distinguish between different fundamental domains in the phase space.

    \subsection{Proposed Method}

    The previous analysis demonstrates that CCM may misinterpret bidirectional causality as unidirectional causality when the underlying dynamic system exhibits two-fold symmetry under the cyclic group $\mathbb{C}_{2}$. The key to resolving this issue lies in reconstructing a one-to-one homeomorphism between two shadow manifolds, enabling CCM to provide correct causality. An intuitive approach is to recover the rotational symmetry of $\mathcal{M}_{z}$, modded by the delay coordinate map, thereby restoring the one-to-one mapping. Existing methods achieve this by incorporating symmetry information from another variable, for example, using $(y,z)$ or $(x,z)$ to reconstruct the Lorenz63 system. However, when applying CCM to detect causality, it is crucial to avoid introducing information from other variables of shadow manifold reconstruction. In other words, we cannot use the shadow manifold reconstructed from the set $(x,z)$ to detect the causality between \(X\) and \(Z\). We find another way to recover the one-to-one mapping without introducing additional information about another variable and propose the segment convergent cross mapping (sCCM) method to address this issue.

    Instead of recovering the modded-out symmetry to reconstruct the one-to-one mapping, we propose another way by projecting the inversion-symmetric shadow manifold (i.e. $\mathcal{M}_{x}$ of Lorenz63) into its fundamental domain $\mathcal{D}_{P}$ and its reflected counterpart $P\cdot\mathcal{D}_{P}$. This segmentation results in two equivalent time sets:
    \begin{equation}
    \begin{aligned}
        [t_{1}] = \{t\in T| \mathbf{x}_{t} \in D_{R_{z}}\},\quad \quad[t_{2}] = \{t\in T| \mathbf{x}_{t} \in R_{z}\cdot{D_{R_{z}}}\}. 
    \end{aligned}
    \end{equation}
    Moreover, two equivalent sub-manifolds, $\mathcal{M}_{x\vert [t_{1}]}$ and $\mathcal{M}_{x\vert [t_{2}]}$, depending on which fundamental domain the trajectory is in at time $t$. Next, we split the invariant shadow manifold $\mathcal{M}_{z}$ into two sub-manifolds $\mathcal{M}_{z\vert [t_{1}]}$ and $\mathcal{M}_{z\vert [t_{2}]}$, corresponding to the time sets $[t_{1}]$ and $[t_{2}]$. In other words, the segmentation of $\mathcal{M}_{z\vert [t_{i}]}$ depends on the indices of $\mathcal{M}_{x\vert [t_{i}]}$, the final step is applying CCM on those two pairs of sub-manifolds for detecting the causality. Through this segmentation process, we eliminate the influence of symmetry, which makes the reconstruction mapping two-to-one (under the nicest case) during the reconstruction process. Moreover, we preserve the local geometry structure of each sub-manifold during the segmentation process. If the bidirectional causality exists, that is the one-to-one mapping $\Tilde{f}: \mathcal{M}_{x\vert [t_{i}]} \rightarrow \mathcal{M}_{z\vert [t_{i}]}$ exists, then the forecasting skill score $\rho$ obtained on two separate experiments should both high and closed to one. Here, we calculate the average as the final score as the final score.

    In a practical scenario, without any analytic information about the dynamic system, the first step is to check the recurrence property of the dynamic system via the recurrence plot or the method shown in \cite{butler2023causal}, which is a necessary condition for the usage of CCM and details are discussed in Appendix \ref{Appendix B}. Next, it is necessary to check the symmetry property of the time series by checking whether the shadow manifolds pose inversion symmetry, which can be achieved by using the indicator shown in \cite{marghoti2024involution}. If one poses inversion symmetry while the other is non-symmetric, our method should be applied under this circumstance to get the correct causality.

\begin{figure}[!h]
        \centering
        \includegraphics[width=14cm]{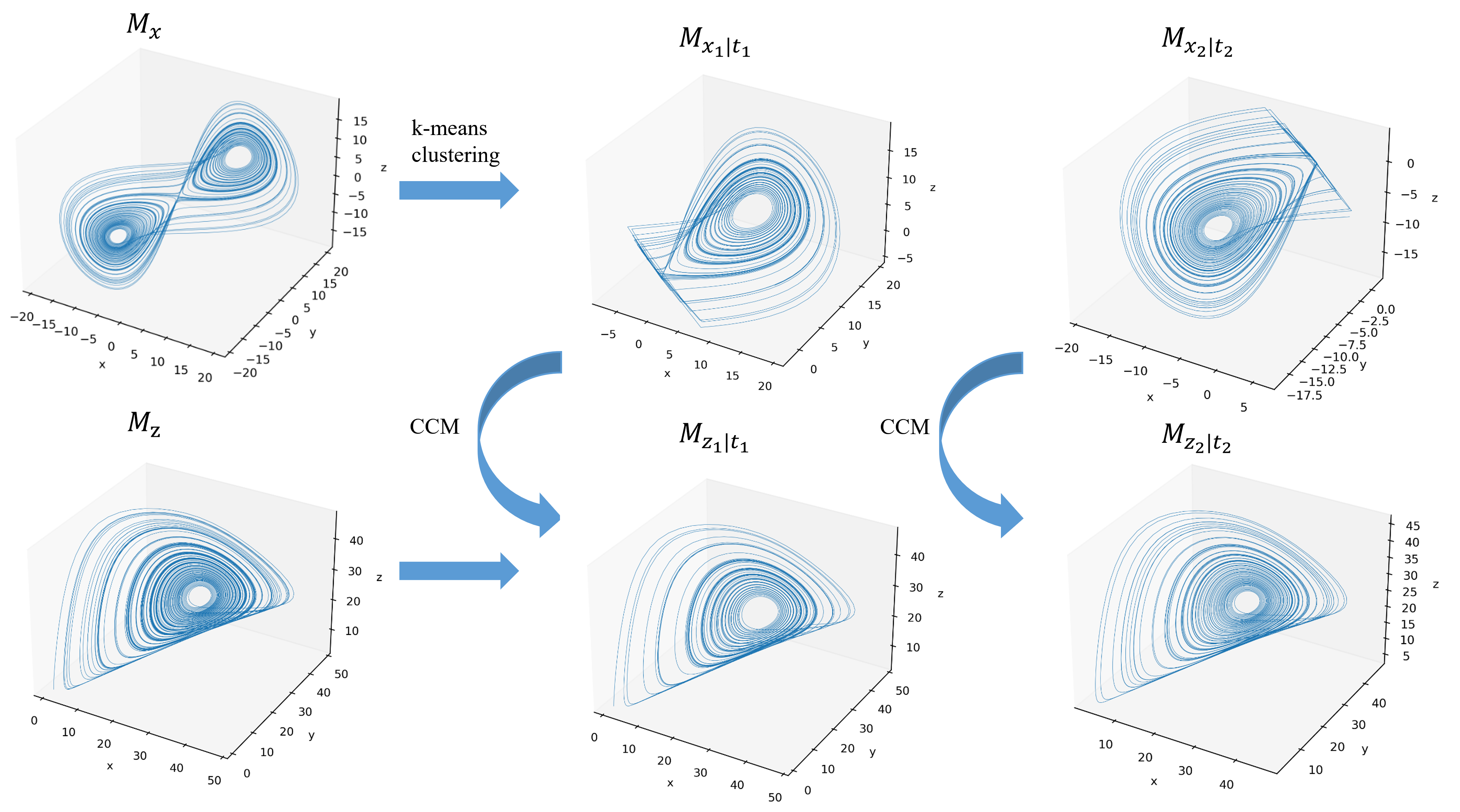}
        \caption{A schematic depicting the details for implementation of sCCM between $\mathcal{M}_{x}$ and $\mathcal{M}_{z}$ for the Lorenz63.}
        \label{Segment Lorenz63}
    \end{figure}

    Here, we use the k-means clustering method as the implementation for the segmentation since k-means performs well on datasets with symmetric properties \cite{ali2022model}. For the Lorenz63 system, the pipeline is illustrated in Fig. \ref{Segment Lorenz63}. Specifically, we first segment the symmetric shadow manifold $\mathcal{M}_{x}$ into two sub-manifolds, $\mathcal{M}_{x\vert [t_{1}]}$ and $\mathcal{M}_{x\vert [t_{2}]}$, using k-means clustering methods (k=2). Subsequently, the indices derived from $\mathcal{M}_{x\vert [t_{i}]}$ are applied to segment $\mathcal{M}_{z}$ into $\mathcal{M}_{z\vert [t_{1}]}$ and $\mathcal{M}_{z\vert [t_{2}]}$.

    \begin{algorithm}[!h]
    \caption{Segment Convergent Cross Mapping for two-fold symmetric system under $\mathbb{C}_{2}$}
    \label{Segment}
    \renewcommand{\algorithmicrequire}{\textbf{Input:}}
    \renewcommand{\algorithmicensure}{\textbf{Output:}}
    \begin{algorithmic}[1]
    \REQUIRE Give two data sets $\{h_{x}\}$ and $\{h_{z}\}$ generated from the coordinate projections, and embedding parameters $\tau$ and $n$  
    \ENSURE Causality between $X$ and $Z$    
    \STATE  Test the recurrence property for the dynamic system. 
    \STATE  Use delay-coordinate mappings $F_{x, \tau, 3}$ and $F_{z, \tau, 3}$ to embed those two data sets into two shadow manifolds $\mathcal{M}_{x}$, and $\mathcal{M}_{z}$.
    \IF {One is inversion symmetric $\mathcal{M}_{x}$ and the other is non-symmetric $\mathcal{M}_{z}$}
    \STATE Use k-means clustering method to partition the symmetric shadow manifold $\mathcal{M}_{x}$ into two parts, $\mathcal{M}_{x\vert [t_{1}]}$ and $\mathcal{M}_{x\vert [t_{2}]}$. 
    \STATE Segment $\mathcal{M}_{z}$ into $\mathcal{M}_{z\vert [t_{1}]}$ and $\mathcal{M}_{z\vert [t_{2}]}$ according to the time index derived from $\mathcal{M}_{x\vert [t_{i}]}$. 
    \STATE Do CCM for $(\mathcal{M}_{x\vert [t_{1}]},  \mathcal{M}_{z\vert [t_{1}]})$; $(\mathcal{M}_{x\vert [t_{2}]}, \mathcal{M}_{z\vert [t_{2}]})$ and obtain $\rho_{x_{1}z_{1}}$, $\rho_{z_{1}x_{1}}$, $\rho_{x_{2}z_{2}}$, $\rho_{z_{2}x_{2}}$
    \STATE Combine the results as $\rho_{xz} = \frac{\rho_{x_{1}z_{1}}+\rho_{x_{2}z_{2}}}{2}$ and $\rho_{zx} = \frac{\rho_{z_{1}x_{1}}+\rho_{z_{2}x_{2}}}{2}$
    \ELSE
    \STATE Do CCM on $\mathcal{M}_{x}$ and $\mathcal{M}_{z}$.
    \ENDIF
    \end{algorithmic}
    \end{algorithm}
     
    Another point to note is that the shadow manifold is when the length of the signal is limited, the number of trajectories lying in the fundamental domain $\mathcal{D}_{P}$ and its reflection counterpart $P\cdot \mathcal{D}_{P}$ of the inversion symmetric shadow manifold may differ; that is, one domain may be denser than the other. This can be clearly observed in Fig. \ref{Segment Lorenz63}, where $\mathcal{M}_{x\vert [t_{2}]}$ looks denser than $\mathcal{M}_{x\vert [t_{1}]}$. This difference is due to the asymmetric periodic motions of the Lorenz63 system; however, as the number of observations increases, both $\mathcal{M}_{x\vert [t_{i}]}$ and $\mathcal{M}_{z\vert [t_{i}]}$ would become denser to the same extent.

    For attractors in higher-dimensional space than three, the complexity far exceeds that of three-dimensional attractors, and the observability dramatically affects the effectiveness of CCM since the quality of reconstructed shadow manifolds faces substantial challenges due to the intricate complicated dynamics. For instance, the system is often unobservable from a single variable, and multiple variables are usually required to obtain an effective embedding \cite{gonzalez2020assessing}. More details are discussed in Appendix \ref{Appendix A}.


    \section{Experiments}\label{SectionFive}
    
    In this section, we test the efficiency of our method through several two-fold rotation symmetric benchmark models under the order-two cyclic group $\mathbb{C}_{2} = \{e, R_{z}(\pi)\}$ in $\mathbb{R}^{3}$. Moreover, we also discuss symmetric systems in the higher dimension.

    \subsection{Case 1: Rotation symmetric systems under order-two group $\mathbb{C}_{2}$ in $\mathbb{R}^{3}$}

    We implement our method on several $\mathbb{C}_{2}$ symmetric systems in $\mathbb{R}^{3}$, including the counterexamples mentioned in Section \ref{SectionThree} and causality results between \(X\) and \(Z\) are shown in Table \ref{Details of Lorenz-like system}. Similarly, the correct relationship between \(Y\) and \(Z\) can be obtained using the same steps. The Pearson correlation coefficients $\rho_{x_{1}z_{1}}$, $\rho_{x_{2}z_{2}}$, $\rho_{z_{1}x_{1}}$, $\rho_{z_{2}x_{2}}$ all converge to nonzero values close to 1 as the length of the library increases, indicating a bidirectional causal relationship between \(X\) and \(Z\).  
     We can observe that CCM misleads the bidirectional causality between two signals as unidirectional causality due to the influence of modded-out symmetry, while sCCM provides the correct causality.

    \begin{table}[h]
        \centering
        \begin{tabular}{c|c|c}
        \hline
        Dynamic systems & CCM & Segment CCM \\
        \hline
        Lorenz63 & $\rho_{xz} = 0.471$, $\rho_{zx} = 0.995$ & $\rho_{xz} = 0.992$, $\rho_{zx} = 0.997$ \\
        \hline
        Chen \& Ueta & $\rho_{xz} = 0.244$, $\rho_{zx} = 0.997$ & $\rho_{xz} = 0.965$, $\rho_{zx} = 0.997$ \\
        \hline
        Burke \& Shaw & $\rho_{xz} = 0.283$, $\rho_{zx} = 0.999$ & $\rho_{xz} = 0.989$, $\rho_{zx} = 0.999$ \\
        \hline
        Three scroll chaotic attractor & $\rho_{xz} = 0.131$, $\rho_{zx} = 0.928$ & $\rho_{xz} = 0.869$, $\rho_{zx} = 0.926$ \\
        \hline
        Wang \cite{wang1992controlling} & $\rho_{xz} = 0.013$, $\rho_{zx} = 0.999$ & $\rho_{xz} = 0.994$, $\rho_{zx} = 0.999$ \\
        \hline
        Shimizu \& Morioka \cite{shimizu1980bifurcation} & $\rho_{xz} = 0.535$, $\rho_{zx} = 0.983$ & $\rho_{xz} = 0.965$, $\rho_{zx} = 0.982$ \\
        \hline
        Rucklidge \cite{rucklidge1992chaos} & $\rho_{xz} = 0.392$, $\rho_{zx} = 0.973$ & $\rho_{xz} = 0.939$, $\rho_{zx} = 0.973$ \\
        \hline
        Sprott B \cite{sprott1994some} & $\rho_{xz} = 0.034$, $\rho_{zx} = 0.999$ & $\rho_{xz} = 0.995$, $\rho_{zx} = 0.999$ \\
        \hline
        Sprott C \cite{sprott1994some} & $\rho_{xz} = 0.326$, $\rho_{zx} = 0.995$ & $\rho_{xz} = 0.970$, $\rho_{zx} = 0.996$ \\
        \hline
        Rikitake \cite{rikitake1958oscillations} & $\rho_{xz} = 0.582$, $\rho_{zx} = 0.991$ & $\rho_{xz} = 0.988$, $\rho_{zx} = 0.987$ \\
        \hline
        Liu \& Yang \cite{lu2004new} & $\rho_{xz} = 0.358$, $\rho_{zx} = 0.989$ & $\rho_{xz} = 0.983$, $\rho_{zx} = 0.967$ \\
        \hline
        Chongxin \textit{et al.} \cite{chongxin2006new} & $\rho_{xz} = 0.185$, $\rho_{zx} = 0.998$ & $\rho_{xz} = 0.968$, $\rho_{zx} = 0.999$ \\
        \hline
        L\"u, Chen \& Cheng \cite{lu2002new} & $\rho_{xz} = 0.035$, $\rho_{zx} = 0.986$ & $\rho_{xz} = 0.946$, $\rho_{zx} = 0.986$ \\
        \hline
        \end{tabular}
        \caption{Results for a range of $\mathbb{C}_{2}$ symmetric systems in $\mathbb{R}^{3}$.}
        \label{Details of Lorenz-like system}
    \end{table}

    Moreover, we test the robustness of sCCM by varying noise levels in the Lorenz63 system. Gaussian noise $n\sim\mathcal{N}(\mu, \sigma)$ is added to the observations. In our tests, $\mu$ is set to zero, and we explore different noise levels with standard deviations $\sigma$ at 0.1, 0.5, and 1. The trajectories are generated as described in Section \ref{SectionThree}, ensuring consistency in other experiments. Results are shown in Table \ref{Noise}.      
    \begin{table}[ht]
        \centering
        \begin{tabular}{c|c|c}
        \hline
        noise level (value of $\sigma$) & CCM & Segment CCM \\
        \hline
        \multirow{2}{*}{\centering 0.1} & $\rho_{xz} = 0.314, \rho_{zx} = 0.989$ & $\rho_{xz} = 0.989, \rho_{zx} = 0.991$ \\
                                        & $\rho_{yz} = 0.315, \rho_{zy} = 0.966$ & $\rho_{yz} = 0.992, \rho_{zy} = 0.969$ \\
        \hline
        \multirow{2}{*}{\centering 0.5} & $\rho_{xz} = 0.153, \rho_{zx} = 0.969$ & $\rho_{xz} = 0.969, \rho_{zx} = 0.969$ \\
                                        & $\rho_{yz} = 0.179, \rho_{zy} = 0.919$ & $\rho_{yz} = 0.975, \rho_{zy} = 0.926$ \\
        \hline
        \multirow{2}{*}{\centering 1}   & $\rho_{xz} = 0.115, \rho_{zx} = 0.920$ & $\rho_{xz} = 0.934, \rho_{zx} = 0.914$ \\
                                        & $\rho_{yz} = 0.121, \rho_{zy} = 0.861$ & $\rho_{yz} = 0.943, \rho_{zy} = 0.867$ \\
        \hline
        \end{tabular}
        \caption{Performance of Segment CCM under different noise levels.}
        \label{Noise}
    \end{table}
    The findings indicate that our method maintains considerable effectiveness even under moderate noise levels, underscoring its robustness. This conclusion is consistent with the other examples and scenarios in this study section.

    \subsection{Case 2: Rotation Symmetric Systems under order-two group $\mathbb{C}_{2}$ in High Dimension}
    
    In this part, we apply our method to three high-dimensional, two-fold rotation symmetric systems, including two four-dimensional and one five-dimensional system. These high-dimensional dynamic systems are visualized by projecting them onto lower-dimensional planes. Observability is critical in determining whether a shadow manifold can be successfully reconstructed in such cases. Consequently, poor observability can lead to incorrect results from CCM.

\begin{figure}[h]
        \centering
        \includegraphics[width = 15cm]{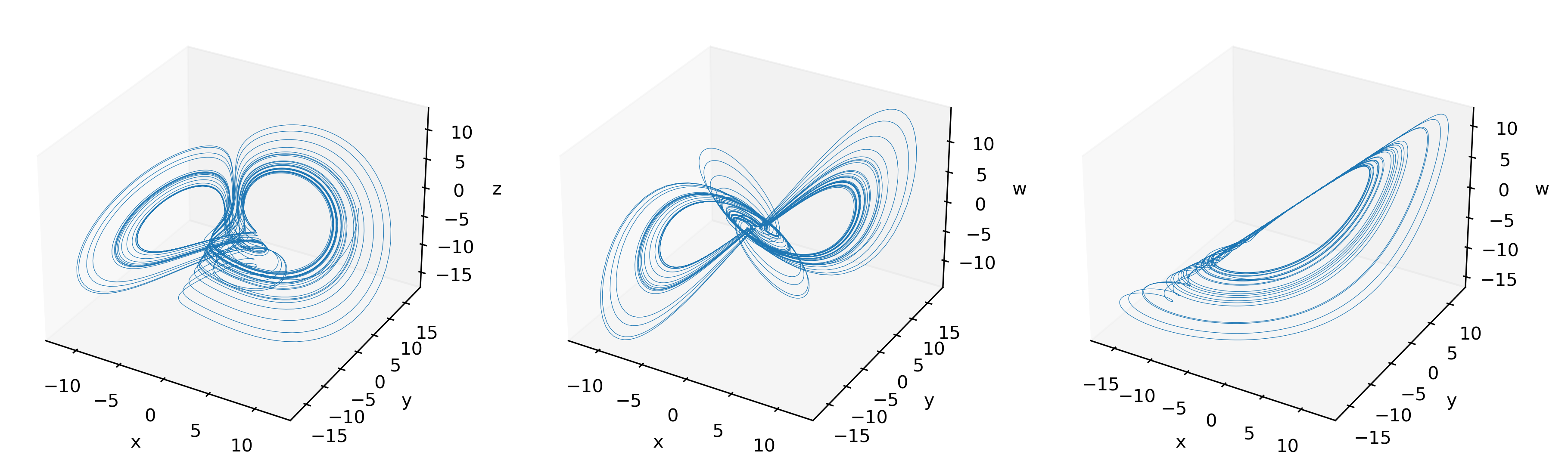}
        \caption{From left to right: Two-wing butterfly chaotic attractors in the $(x,y,z)$ plane, the four-wing butterfly chaotic attractors in the $(x,y,w)$ plane, and the $(x,y,w)$ projection of the reconstructed shadow manifold $\mathcal{M}_{z}$, where $\tau = 50$, $n = 4$.}
        \label{4D Figure}
    \end{figure}

    The first case is a four-dimensional dissipative Lorenz-like dynamic system, representing a self-excited oscillatory modular circuit with the following mathematical formulation \cite{huang2019new}:
    \begin{equation}
        \begin{aligned}
            & \frac{dx}{dt} = a(y-x)\\
            & \frac{dy}{dt} = xz + w\\
            & \frac{dz}{dt} = b-xy\\
            & \frac{dw}{dt} = yz-cw\\
        \end{aligned}\label{4D}
    \end{equation}
    where parameters are set as $a =6, b=11, c=5$, with the initial value $(x,y,z,w) = (10,10,0,0)$. We simulate this system using the fourth-order Runge-Kutta method over the temporal domain $t=[0,100.0]$ with a time step of $\Delta t = 0.01$. This system is rotation symmetric under the order-two group $\mathbb{C}_{2} = \{e, R_{z}(\pi)\}$, where $R_{z}:(x, y, z, w) \rightarrow (-x, -y, z,  -w)$, producing two distinct types of attractors: four-wing and two-wing attractors which are shown in Fig. \ref{4D Figure}.

    From equation (\ref{4D}), we infer bidirectional causal links between variables: $X\Leftrightarrow Z$, $Y\Leftrightarrow Z$, and $W\Leftrightarrow Z$. However, CCM results show unidirectional causal links $Z\Rightarrow X$, $Z\Rightarrow Y$, and $Z\Rightarrow W$, because of the nonsymmetric shadow manifold $\mathcal{M}_{z}$. By applying sCCM, we successfully recover the correct causalities, and results are shown in Table \ref{High dimensional table}.

    Another four-dimensional Lorenz-like hyperchaotic dynamic system is given by \cite{yang2009hyperchaotic}:
    \begin{equation}
        \begin{aligned}
            & \frac{dx}{dt} = 10(y-x)\\
            & \frac{dy}{dt} = 28x - y - xz + w\\
            & \frac{dz}{dt} = -\frac{8}{3}z + xy\\
            & \frac{dw}{dt} = -k_{1}x -k_{2}y\\
        \end{aligned}\label{4D Lorenz}
    \end{equation}
    where $k_{1} = -9.3$, and $k_{2} = -5$. This dynamic system is constructed by adding a linear controller to the second equation of the Lorenz63 system. This dynamic system is rotational symmetric under the order-two group $\mathbb{C}_{2} = \{e, R_{z}\}$, where $R_{z}: (x, y, z, w) \Rightarrow (-x, -y, z, w)$. We simulate this dynamic system in the temporal domain $t = [0, 50]$ with $\Delta t = 0.005$, and parameters for the embedding are $\tau = 15$, $n = 4$. Similarly, all shadow manifolds except $\mathcal{M}_{z}$ are symmetric. Consequently, CCM incorrectly identifies bidirectional causal links $X \Leftrightarrow Z$, $Y \Leftrightarrow Z$, and $W \Leftrightarrow Z$ as unidirectional. True causal links are correctly detected by applying sCCM, and results are shown in Table \ref{High dimensional table}.

\begin{figure}[h]
        \centering
        \includegraphics[width=14cm]{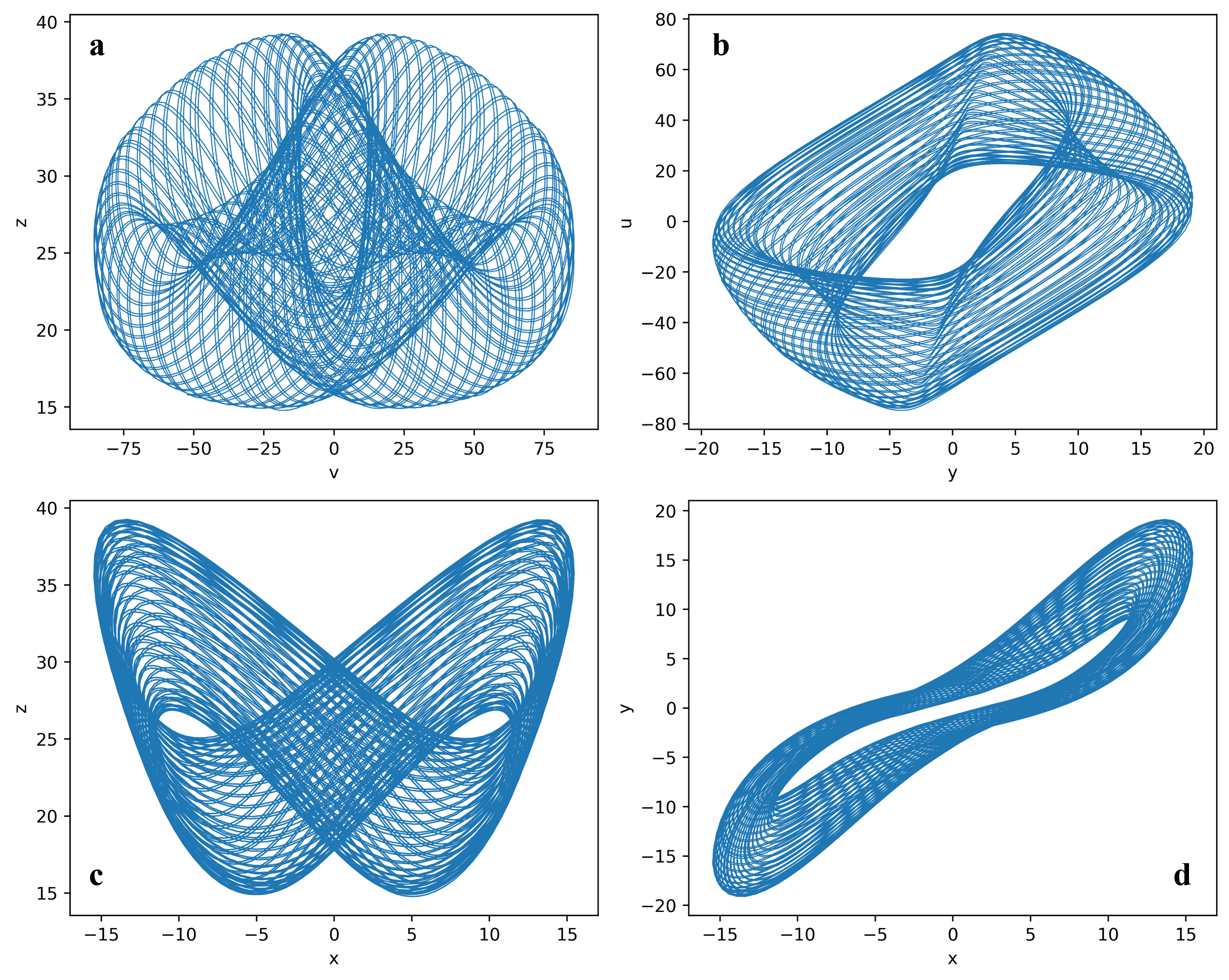}
        \caption{Projections of the attractor generated by equation (\ref{5D formulation}) a. (v,z) plane projection. b. (y,u) plane projection. c. (x,z) plane projection. d. (x,y) plane projection.}
        \label{5D}
    \end{figure}

    The last case is a five-dimensional hyperchaotic system. Similarly, this system is obtained by adding a nonlinear quadratic controller to the first equation and a linear controller to the second equation of the modified Lorenz63 system. Moreover, this system can be realized via an electronic system \cite{hu2009generating}. It is mathematical formulation is:
    \begin{equation}
        \begin{aligned}
            & \frac{dx}{dt} = 10(y-x) + u\\
            & \frac{dy}{dt} = 28x - y - xz - v\\
            & \frac{dz}{dt} = -\frac{8}{3}z + xy\\
            & \frac{du}{dt} = -xz + k_{1}u\\ 
            & \frac{dv}{dt} = k_{2}y\\
        \end{aligned}\label{5D formulation}
    \end{equation}
    where $k_{1} = 1$ and $k_{2} = 30$. This system is simulated in the temporal domain $t = [0, 100]$ with $\Delta t = 0.01$, and embedding parameters are $\tau = 60$, $ n = 5$. This hyperchaotic system is rotation symmetric under the order-two group $\mathbb{C}_{2} = \{e, R_{z}(\pi)\}$, where $R_{z}: (x, y, z, u, w) \Rightarrow (-x, -y, z, -u, -w)$ as shown in Fig. \ref{5D}, which exhibits a genus-1 attractor in the $(v,z)$ plane.  
    
    Similarly, the nonsymmetric shadow manifold $\mathcal{M}_{z}$ causes CCM to mislead bidirectional causal links $X\Leftrightarrow Z$, $Y\Leftrightarrow Z$, and $W\Leftrightarrow Z$ as $Z\Rightarrow X$, $Z\Rightarrow Y$, and $Z\Rightarrow W$. However, true causal links are correctly identified using our framework. Results are shown in Table \ref{High dimensional table}.     

    \begin{table}[ht]
    \centering
    \begin{tabular}{c|c|c}
    Equations & CCM & Segment CCM \\
    \hline
    \multirow{3}{*}{4D dissipative system (\ref{4D})}                                    & $Z \Rightarrow X$ $(\rho_{xz} = 0.547, \rho_{zx} = 0.994)$ & $X \Leftrightarrow Z$ $(\rho_{xz} = 0.967, \rho_{zx} = 0.973)$\\ 
                       & $Z \Rightarrow Y$ $(\rho_{yz} = 0.474, \rho_{zy} = 0.990)$ 
    & $Y \Leftrightarrow Z$ $(\rho_{yz} = 0.971, \rho_{zy} = 0.931)$ \\  
                       & $Z \Rightarrow W$ $(\rho_{wz} = 0.273, \rho_{zw} = 0.978)$ 
    & $W \Leftrightarrow Z$ $(\rho_{wz} = 0.950, \rho_{zw} = 0.904)$ \\
    \hline
    \multirow{3}{*}{4D Lorenz-like system (\ref{4D Lorenz})}                            & $Z \Rightarrow X$ $(\rho_{xz} = 0.597, \rho_{zx} = 0.988)$ & $X \Leftrightarrow Z$ $(\rho_{xz} = 0.973, \rho_{zx} = 0.964)$\\ 
                       & $Z \Rightarrow Y$ $(\rho_{yz} = 0.583, \rho_{zy} = 0.962)$
    & $Y \Leftrightarrow Z$ $(\rho_{yz} = 0.966, \rho_{zy} = 0.905)$ \\  
                       & $Z \Rightarrow W$ $(\rho_{wz} = 0.609, \rho_{zw} = 0.954)$
    & $W \Leftrightarrow Z$ $(\rho_{wz} = 0.912, \rho_{zw} = 0.861)$ \\
    \hline
    \multirow{4}{*}{5D Lorenz-like system (\ref{5D formulation})}                       & $Z \Rightarrow X$ $(\rho_{xz} = 0.151, \rho_{zx} = 0.999)$ & $X \Leftrightarrow Z$ $(\rho_{xz} = 0.994, \rho_{zx} = 0.999)$ \\  
                       & $Z \Rightarrow Y$ $(\rho_{yz} = 0.194, \rho_{zy} = 0.999)$ & $Y \Leftrightarrow Z$ $(\rho_{yz} = 0.985, \rho_{zy} = 0.986)$ \\  
                       & $Z \Rightarrow U$ $(\rho_{uz} = 0.098, \rho_{zu} = 0.999)$ & $U \Leftrightarrow Z$ $(\rho_{uz} = 0.991, \rho_{zu} = 0.999)$ \\
                       & $Z \Rightarrow W$ $(\rho_{wz} = 0.105, \rho_{zw} = 0.999)$ & $W \Leftrightarrow Z$ $(\rho_{wz} = 0.993, \rho_{zw} = 0.999)$\\
    \hline
    \end{tabular}
    \caption{Performance of Segment CCM for high dimensional rotational symmetric systems under $\mathbb{C}_{2}$.}
    \label{High dimensional table}
    \end{table}
    
    \subsection{Further discussion: Rotation Symmetric Systems under order-$n$ group $\mathbb{C}_{n}$}

    As we mentioned in Section \ref{SectionFour}, if the shadow manifold reconstructed using delay-coordinate mapping exhibits symmetry, it can, at most, preserve two-fold inversion symmetry, which implies that the reconstructed shadow manifold cannot maintain the original attractor's k-fold symmetry if $k > 2$. In this case, an obvious question arises: when the original attractor exhibits k-fold symmetry (k>2), whether the reconstructed symmetric shadow manifold only preserves two-fold inversion symmetry may affect the result of CCM. We found that in some cases, even if the same order of symmetry as the original attractor cannot be preserved, the delay-coordinate mapping used to reconstruct the shadow manifold is a one-to-one mapping. Thus, it does not affect the results of CCM. To further illustrate this point, we use the four-fold Burke \& Shaw system \cite{gilmore2007symmetry}, a rotational symmetric system under $\mathbb{C}_{4}$ as an example.

\begin{figure}[h]
        \centering
        \includegraphics[width=16cm]{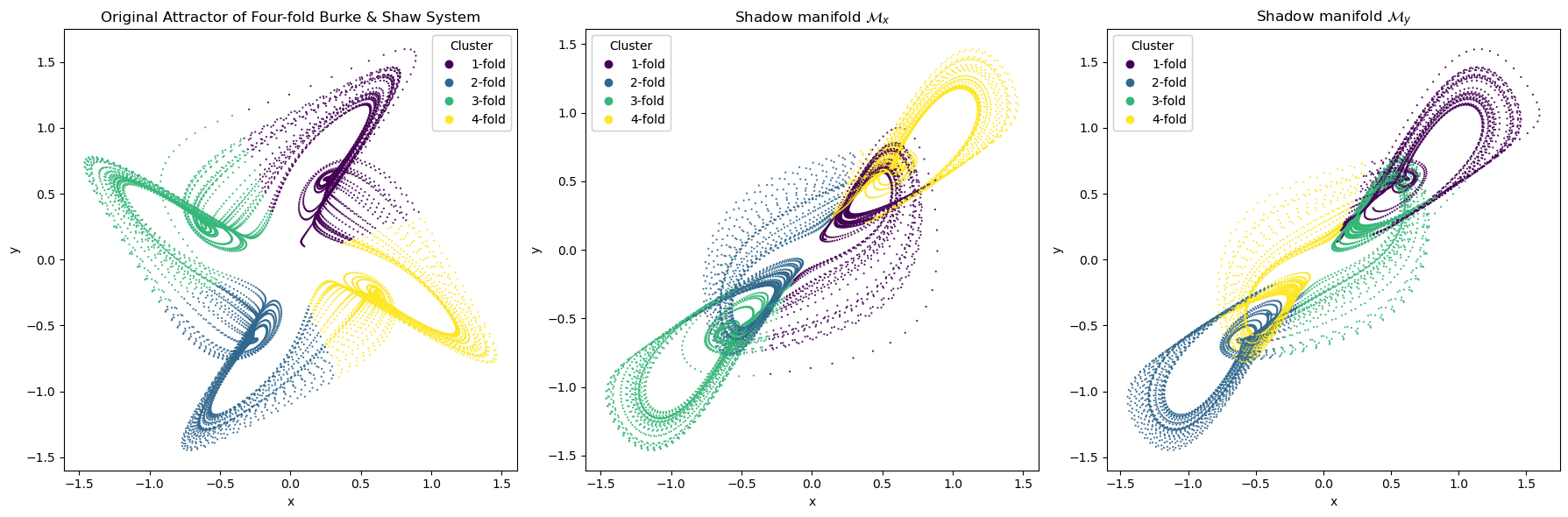}
        \caption{Attractors of equation (\ref{Four-fold Burke and Shaw}). a. $(x,y)$ plane projection of the four-fold Burke $\&$ Shaw Attractor. b. $(x,y)$ plane projection of the shadow manifold $\mathcal{M}_{x}$. c. $(x,y)$ plane projection of the shadow manifold $\mathcal{M}_{y}$. }
        \label{Four-fold}
    \end{figure}

    The mathematical formulation of the four-fold Burke \& Shaw system is:
    \begin{equation}\label{Four-fold Burke and Shaw}
    \begin{aligned}
        \dot{x} &= -\frac{(S + 1)x}{4} - \frac{S(1 - z)y}{4} + \frac{u_3(1 - S) - v_3S(1 + z)}{4R^2}, \\
        \dot{y} &= \frac{S(1 - z)x}{4} - \frac{(S + 1)y}{4} - \frac{v_3(1 - S) + u_3S(1 + z)}{4R^2}, \\
        \dot{z} &= \mathcal{V} + \frac{S}{2}v_4.
    \end{aligned}
    \end{equation}
    where $R^{2} = x^{2} + y^{2}$, $u_{3} = x^3 -3xy^{2}$, $v_{3} = 3x^{2}y - y^{3}$, $v_{4} = 4x^{3}y - 4xy^{3}$, $S = 10$ and $\mathcal{V} = 4.271$. This system is simulated in the temporal domain $t = [0, 100]$, $\Delta t = 0.01$ with initial point $\mathbf{x}_{0} = [0.1, 0.1, 0.1]$, and embedding parameters are $\tau = 10$ and $n = 3$. The image of this four-fold symmetric system is shown in Fig. \ref{Four-fold}. To better illustrate that the delay-coordinate mapping is one-to-one, we assign different colors to the trajectories on each cover of the original four-fold attractor. The points on the shadow manifold are then passively assigned colors based on the corresponding points on the original attractor. From the two-fold symmetric shadow manifold $\mathcal{M}_{x}$, the blue and purple cover degenerate as two symmetric bonuses linking the other two covers. The same result can also be obtained for $\mathcal{M}_{y}$. Although the original four-fold symmetry degenerates to the two-fold symmetry because of the delay-coordinate mappings $F_{x, \tau. 3}$ and $F_{z, \tau, 3}$, these mappings are still one-to-one. Thus, the induced homeomorphism $f$ between those two shadow manifolds can also be induced, which may not affect the causality detection of CCM between those two signals. Detailed results are shown in Table \ref{N-fold Table}. However, under this circumstance, sCCM does not perform well since it can not reveal the one-to-one mapping through the segmentation process through k-means clustering because of the degeneration of covers. Further discussion is needed to find a refinement segmentation framework. 

    \begin{table}[h]
        \centering
        \begin{tabular}{c|c}
        \hline
        Equations & \multicolumn{1}{c}{CCM} \\
        \hline
        Four-fold Burke $\&$ Shaw System
            & \multicolumn{1}{c}{$X \Leftrightarrow Y \ (\rho_{xy} = 0.981, \ \rho_{yx} = 0.982)$} \\
        \hline
        \end{tabular}
        \caption{Detailed results for the causality test for the circuit case.}
        \label{N-fold Table}
    \end{table}

    \section{Conclusion}\label{Conclusion}
    
    This paper explores the complexities of detecting causal relations in nonlinear dynamic systems using CCM, explicitly focusing on the challenges of rotation symmetric systems. The result reveals that shadow manifolds reconstructed by the differential/delay-coordinate mapping from the invariant variable, lacking symmetry information, can potentially mislead the causality detection in CCM since the delay-coordinate mapping is no longer a one-to-one mapping. For systems that exhibit two-fold rotational symmetry, we proposed a new framework based on the k-means clustering method, partitioning the symmetric shadow manifolds into two sub-manifolds, thus transforming the two-to-one mapping into two one-to-one. We evaluate the efficiency of our method on various low and high-dimensional two-fold rotation symmetric dynamic systems, as well as a real-world application involving a nonlinear circuit, to verify its accuracy.

    Several areas warrant further exploration. First, the result of CCM is highly dependent on the performance of the delay-coordinate mapping. As we discussed, for high-dimensional dynamic systems, one-dimensional observations generally provide poor observability of the original attractor, necessitating the combination of observations from multiple variables for sufficient embedding. However, this approach may introduce confounding information, complicating the detection of causal relationships. Second, our method may not perform well for a k-fold ($k>2$) symmetric system, and a further refined segmentation framework should be proposed. 
    
    \section*{Acknowledgments}
    
    The author acknowledge that no external funding or support was received for this research.

    \bibliographystyle{elsarticle-num}
    \bibliography{main}

    \section*{Appendix A: Observability of nonlinear dynamic system}\label{Appendix A}

    According to Taken's embedding theorem, using a delay-coordinate mapping, an $n$-dimensional dynamical system can be reconstructed from a one-dimensional observation. This reconstruction process is akin to splitting $n$-dimensional linearly independent information. However, the effectiveness of this approach diminishes as the dimension of the dynamical system increases. To illustrate these limitations, we consider the nine-dimensional Lorenz system \cite{reiterer1998nine} as an example, demonstrating the inefficient performance of the delay-coordinate mapping in high-dimensional settings.

\begin{figure}[h]
        \centering
        \includegraphics[width=14cm]{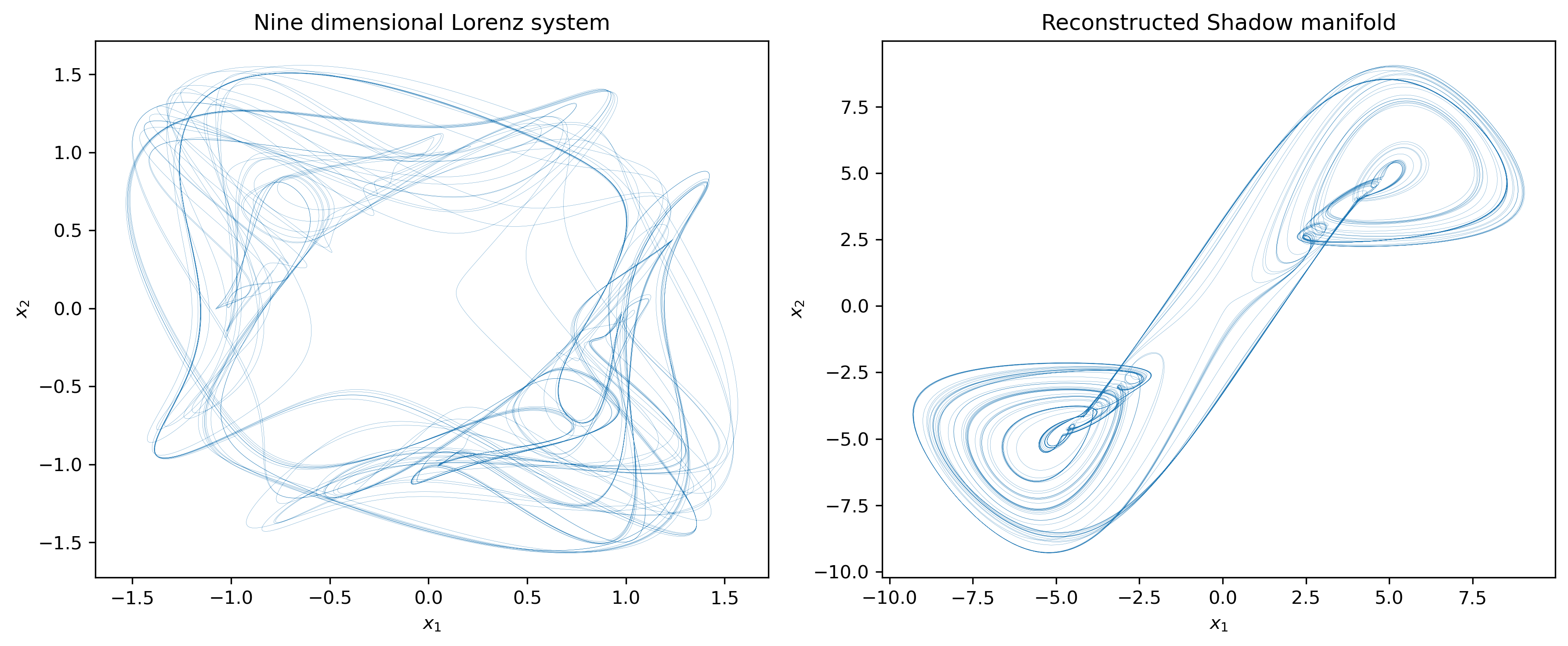}
        \caption{\textbf{Left: Chaotic attractor produced by nine-dimensional Lorenz system. Right: Reconstructed shadow manifold by delay-coordinate mapping. }}
        \label{9D Lorenz}
    \end{figure}

    The mathematical formulation of the nine-dimensional Lorenz system is given by: 
    \begin{equation}
    \begin{aligned}
        \dot{x}_1 &= -\sigma (b_1 x_1 + b_2 x_7) + x_4 (b_4 x_4 - x_2) + b_3 x_3 x_5, \\
        \dot{x}_2 &= -\sigma x_2 + x_1 x_4 - x_2 x_5 + x_4 x_5 - \frac{2}{\sigma} x_9, \\
        \dot{x}_3 &= \sigma (b_2 x_8 - b_1 x_3) + x_2 x_4 - b_4 x_2^2 - \frac{b_3 x_1 x_5}{\sigma}, \\
        \dot{x}_4 &= -\sigma x_4 - x_2 x_3 - x_2 x_5 + x_4 x_5 + \frac{x_9}{2}, \\
        \dot{x}_5 &= -\sigma b_5 x_5 + \frac{x_2^2}{2} - \frac{x_4^2}{2}, \\
        \dot{x}_6 &= -b_6 x_6 + x_2 x_9 - x_4 x_9, \\
        \dot{x}_7 &= -b_1 x_7 - R x_1 + 2 x_5 x_8 - x_4 x_9, \\
        \dot{x}_8 &= -b_1 x_8 + R x_3 - 2 x_5 x_7 + x_2 x_9, \\
        \dot{x}_9 &= -x_9 + (R + 2 x_6)(x_4 - x_2) + x_4 x_7 - x_2 x_8,
    \end{aligned}
    \end{equation}
    \noindent where $b_{1} = \frac{5}{1.5}$, $b_{2} = 0.6$, $b_{3} = 1.2$, $b_{4} = 0.2$, $b_{5} = \frac{2}{1.5}$, $b_{6} = \frac{4}{1.5}$, and $R = 14.3$. The initial state is $\mathbf{x}_{0} = [0.01, 0, 0.01, 0, 0, 0, 0, 0, 0.01]$, with a time step of $\Delta t = 0.02$ and the simulation time domain of $[0, 800]$. The measurement function is $h = x_{9}$, with the lag value $\tau = 12$ and an embedding dimension $n = 9$. The original attractor and the reconstructed shadow manifold are shown in Fig. \ref{9D Lorenz}.
    \noindent It seems that the reconstructed shadow manifold exhibits a butterfly attractor similar to the Lorenz63 system rather than reflecting the more complex dynamic behavior in the nine-dimensional attractor. This discrepancy can be attributed to the observability of the system. Observability refers to the ability to reconstruct the entire state of the system using a subset of measurable variables over a finite period. Consider the following nonlinear system:
    \begin{equation}\label{nonlinear}
    \begin{aligned}
             \Dot{\mathbf{x}} &= \mathbf{f}(\mathbf{x}) \\
             s(t) &= h(\mathbf{x})
    \end{aligned}
    \end{equation} 
    where $\mathbf{x}\in \mathbb{R}^{m}$ represents the $m$-dimensional state vector, $h:\mathbb{R}^{m}\rightarrow \mathbb{R}$ is the measurement function, $s(t)\in \mathbb{R}$ is the measurement, and $\mathbf{f}: \mathbb{R}^{m} \rightarrow \mathbb{R}^{m}$ denotes the vector field. Differentiating $s(t)$ yields 
    \begin{equation}
    \begin{aligned}
             \dot{s}(t) = \frac{d}{dt}h(\mathbf{x}) = \frac{\partial h}{\partial \mathbf{x}}\mathbf{\Dot{x}} = \frac{\partial h}{\partial \mathbf{x}}\mathbf{f}(\mathbf{x}) = \mathcal{L}_{f}h(\mathbf{x})  
    \end{aligned}
    \end{equation}
    where $\mathcal{L}_{f}h(\mathbf{x})$ is the Lie derivative of h along the vector field $\mathbf{f}$ and the time derivative of $s$ in Lie derivatives form can be written as $s^{j}=\mathcal{L}_{f}^{i}h(\mathbf{x})$. The $j$-th order Lie derivative is given by
    \begin{equation}
    \begin{aligned}
        \mathcal{L}_{f}^{j}h(\mathbf{x}) = \frac{\partial \mathcal{L}_{f}^{j-1}h(\mathbf{x})}{\partial \mathbf{x}}\cdot \mathbf{f}(\mathbf{x})
    \end{aligned}
    \end{equation}
    where $\mathcal{L}_{f}^{0}h(\mathbf{x}) = h(\mathbf{x})$.

    A dynamic system is said to be \textit{state observable} at time $t_{f}$ if its initial state $\mathbf{x}_{0}$ can be uniquely determined from the measurement vector $s(t)$ over the interval $0\leq t\leq t_{f}$. The observability matrix, the Jacobian matrix of the Lie derivatives of $h(\mathbf{x})$, determines whether a dynamic system is observable. The observability matrix $\mathcal{O} \in \mathbb{R}^{m\times m}$ is expressed as:
    \begin{equation}
    \mathcal{O}(\mathbf{x})=\left[ \begin{array}{c}
    \frac{\partial \mathcal{L}_{f}^{0}h(\mathbf{x})}{\partial \mathbf{x}} \\
    \vdots \\
    \frac{\partial \mathcal{L}_{f}^{m-1}h(\mathbf{x})}{\partial \mathbf{x}} \\
    \end{array} \right].
    \end{equation}
    For the linear system, 
    \begin{equation}
    \begin{aligned}
        \Dot{\mathbf{x}} &= A\mathbf{x} + B\mathbf{u} \\
        \mathbf{s} &= C\mathbf{x},
    \end{aligned}
    \end{equation}
    where $A, B, C$ are constant dynamics matrices, the observability matrix simplifies to:
    \begin{equation}
    \mathcal{O}=\left[ \begin{array}{c}
    C\quad CA\quad CA^{2}\quad \cdots \quad CA^{n-1} \\
    \end{array} \right]^{T}.
    \end{equation}
    
    \textbf{Theorem 1.} The $m$-dimensional dynamic system is said to be \textit{state observable} if and only if the observability matrix has full rank, that is, $rank(\mathcal{O}) = m$.
    
    The above theorem is a generalization of the observability property of both linear and nonlinear systems and has been proved in \cite{hermann1977nonlinear}. It states that if a system is observable for both linear and nonlinear dynamic systems, it is possible to recover every initial condition from the measured series $s(t)$ for $t\leq 0$. In the case of the nine-dimensional Lorenz system, the system is nearly unobservable from a single observation, as the observability matrix is deficient. For this system, at least six variables should be measured to achieve sufficient observability of the high dimensional dynamics \cite{letellier2018symbolic}.

    Thus, for high-dimensional dynamical systems, a single variable is generally insufficient to construct a diffeomorphic shadow manifold, implying that CCM results may be unreliable.

\begin{figure}[h]
        \centering
        \includegraphics[width=14cm]{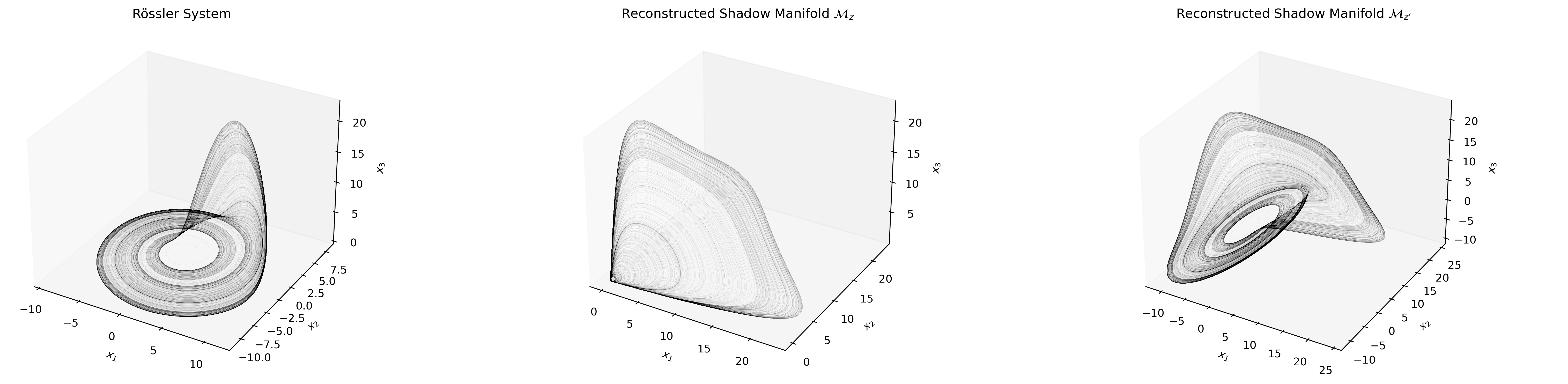}
        \caption{\textbf{Left: Chaotic attractor produced by the R\"{o}ssler system. Middle: Reconstructed shadow manifold $\mathcal{M}_{z}$. Right: Reconstructed shadow manifold $\mathcal{M}_{z'}$}.}
        \label{Rossler}
    \end{figure}

    Moreover, a single variable may still provide poor observability even for low-dimensional dynamic systems ($n\leq 3$). For instance, consider the R\"{o}ssler system, the mathematical formulation is:
    \begin{equation}
    \begin{aligned}
        \dot{x}_1 &= -y-z, \\
        \dot{x}_2 &= x + ay, \\
        \dot{x}_3 &= b + xz -cz, \\
    \end{aligned}
    \end{equation}
    where $a = 0.2$, $b=0.2$, $c=5.7$. The time step is set to $\Delta t = 0.01$, with a total of $400000$ iterations, starting from the initial condition $\mathbf{x}_{0} = [1,1,0]$. The original attractor and the reconstructed shadow manifold are shown in Fig. \ref{Rossler}. 
    Here, we modify the measurement function to $h' = x+z$ and obtain the shadow manifold $\mathcal{M}_{z'}$. The lag value is $\tau =40$ for both shadow manifolds, and the embedding dimension is $n = 3$. It can be observed that the single observation $h=z$ yields poor observability, such that the observability matrix $\mathcal{O}(\mathbf{x})$ is deficient near the original point, leading to significant distortion in the shadow manifold $\mathcal{M}_{z}$. In contrast, the shadow manifold reconstructed from the multivariate measurement function $h'$ provides a better result.

    \section*{Appendix B: The Importance of Recurrence Property: Explanations for a case}\label{Appendix B}

\begin{figure}[ht]
        \centering
        \includegraphics[width=14cm]{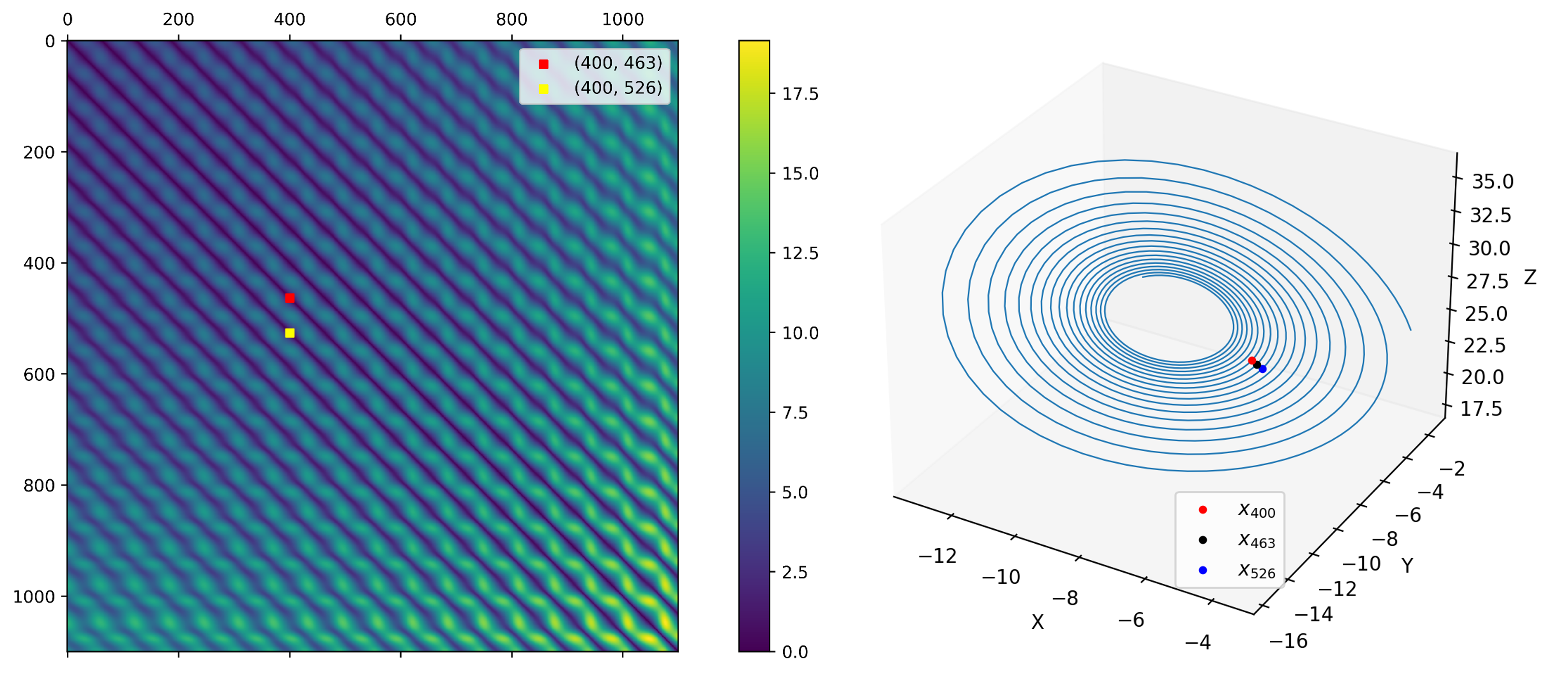}
        \caption{Distance plot using the Euclidean distance and corresponding trajectories of Lorenz63 attractor. This system is simulated as same as Fig. \ref{CCM4Lorenz63}. Left: Distance plot for trajectories lies in $t\in [1, 12]$. Right: Trajectories for the Lorenz63 at $t\in [1, 12]$.}
        \label{RecurrencePlot}
    \end{figure}

    Here, we emphasize that the recurrence property is necessary for CCM. The recurrence property is one of the fundamental properties of the dynamic systems. As stated in the Poincaré recurrence theorem \cite{katok1995introduction}, in any measure-preserving transformation on a dynamic system's attractor, the trajectories will eventually reappear at the neighborhood of the former points in the state space. The recurrences can be captured through distance plots \cite{marwan2007recurrence}. For a series of trajectories $\{\mathbf{x}_{i}\}_{i=1}^{N}$ of a dynamic system in its phase space, the corresponding distance plot is based on the matrix:
    \begin{equation}
        R_{i,j} = \parallel \mathbf{x}_{i}-\mathbf{x}_{j}\parallel,\quad i,j = {1,2,...,N}
    \end{equation}
    where $N$ is the number of considered states $\mathbf{x}_{i}$ and $\parallel\cdot\parallel$ represents a norm (e.g., the Euclidean norm). An example of Lorenz63 system is presented in Fig. \ref{RecurrencePlot},
    it delinates the distances between every pair of points $\mathbf{x_{i}}$ and $\mathbf{x_{j}}$ in the state space. The points $\mathbf{x}_{400}$ and $\mathbf{x}_{526}$ would revisit the position within $\mathbf{x}_{400}$.

    Consider the counterexample mentioned in \cite{bartsev2021imperfection}. This paper discusses that when the time series have certain forms – one of the series is a periodic function of time, and the second series has a time trend - it may mislead the results of CCM because points from the manifold with a periodic shape are close to points from the manifold with a linear shape. However, we point out that the second monotonic increasing time series does not satisfy the recurrence requirement for CCM, which means the shadow manifold for this monotone time series is not an attractor.

    One typical example is:
    \begin{equation}\label{conter}
        \begin{aligned}
            & x(t) = \frac{1}{2000}t\\
            & y(t) = \mathrm{sin}(\frac{\pi}{50}t)
        \end{aligned}
    \end{equation}
    Here $t$ refers to the time. We generate 2000 examples containing 20 complete periods about $y(t)$, and the results of causal links are shown in the following Figure \ref{CounterExample}.
    \begin{figure}[ht]
        \centering
        \includegraphics[width=14cm]{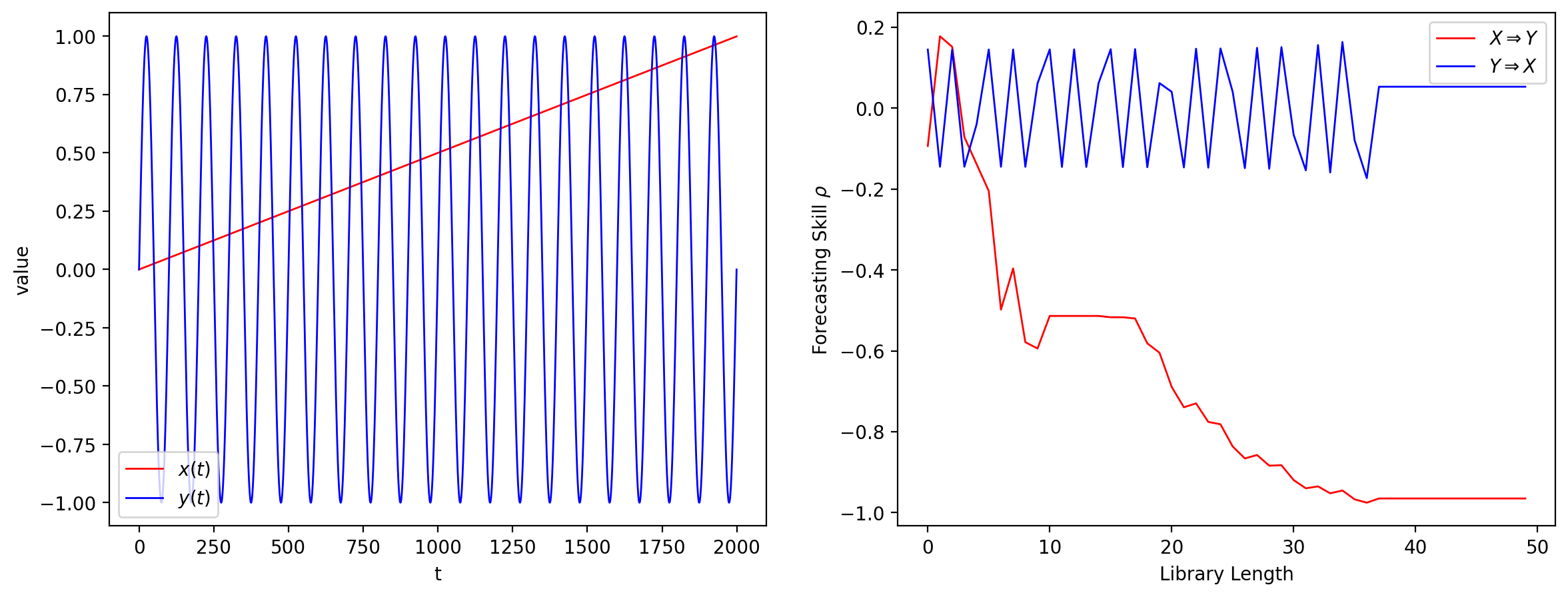}
        \caption{Left: Trajectories of Eq. (\ref{conter}). Right: Causal links obtained by CCM.}
        \label{CounterExample}
    \end{figure}
    The monotone variable causes the periodic one as the Pearson correlation coefficient $rho_ {xy}$ converges to 0.965. In \cite{bartsev2021imperfection}, it is stated that a monotone time series has no causal effect on the periodic time series. However, the mathematical expression of $y(t)$ can be simply solvable by $x(t)$ as we wrote it into the discrete form:
    \begin{equation}
    \begin{aligned}
            x(t+1) &= x(t) + \frac{1}{2000}\\
            y(t+1) &= \mathrm{sin}(40\pi[x(t+1)-\frac{1}{2000}]+\frac{1}{50})\\
    \end{aligned}
    \end{equation}
    which means that the monotone variable causes the other one and corresponds to what we see. However, the inverse relationship also holds:
    \begin{equation}
    \begin{aligned}
            y(t+1) &= y(t) \cos\left( \frac{\pi}{50} \right) + \cos\left( \frac{\pi}{50} t \right) \sin\left( \frac{\pi}{50} \right)\\
            x(t+1) &= \frac{\arcsin(y(t+1)) - \frac{1}{50}}{40\pi} + \frac{1}{2000}\\
    \end{aligned}
    \end{equation}
    From the above analysis, if the nonlinear periodic term can be explicitly solved, as in the case above, it can be expressed in terms of the linear term and vice versa. However, it is important to emphasize that the CCM test becomes nonsensical under these conditions. It does not conform to the requirement of CCM; precisely, with the library length increase, the monotone shadow manifold reconstructed from the linear term would not become denser with the nonlinear one. In this scenario, the prediction ability does not improve as the library length increases, making the CCM result unreliable.

\end{document}